\documentclass[7pt]{article}                    % onecolumn (standard format)
\usepackage{float}
\usepackage{graphicx,epstopdf}
\usepackage{bbm}
\usepackage{amsfonts}
\usepackage{booktabs}
\usepackage{multirow}
\usepackage{siunitx}
\usepackage{mathrsfs}
\usepackage{csquotes}
\usepackage{fancyhdr}
\usepackage{amsmath}
\usepackage{flexisym}
\usepackage{dsfont}
\usepackage{amssymb}
\usepackage{theorem}
\usepackage{graphicx}
\usepackage[dvips]{epsfig}
\usepackage{latexsym}
\usepackage{exscale}
\usepackage[latin1]{inputenc}

\newtheorem{theorem}{Theorem}[section]
\newtheorem{lemma}{Lemma}[section]

\newtheorem{proposition}{Proposition}[section]

\newtheorem{definition}{Definition}[section]

\newtheorem{corollary}{Corollary}[section]

\numberwithin{equation}{section}

\begin{document}

\title{Estimation of the Pointwise H\"older Exponent of Hidden Multifractional Brownian Motion Using Wavelet Coefficients%\thanks{Grants or other notes
}

%\titlerunning{Estimation of Pointwise H\"older Exponent}        % if too long for running head

\author{Sixian Jin\footnote{ Institute of Mathematical Sciences, Claremont Graduate University, USA. E-mail: Sixian.Jin@cgu.edu} \and Qidi Peng\footnote{Corresponding author, Institute of Mathematical Sciences, Claremont Graduate University, USA. E-mail: Qidi.Peng@cgu.edu} \and Henry Schellhorn\footnote{Institute of Mathematical Sciences, Claremont Graduate University, USA. E-mail: Henry.Schellhorn@cgu.edu}}

%\authorrunning{S. Jin, Q. Peng and H. Schellhorn} % if too long for running head

\date{}

\maketitle
\begin{abstract}
We propose a wavelet-based approach to construct consistent estimators of the pointwise H\"older exponent of a multifractional Brownian motion, in the case where this underlying process is not directly observed. The relative merits of our estimator are discussed, and we introduce an application to the problem of estimating the functional parameter of a nonlinear model.\\
\textbf{keywords:}{ Pointwise H\"older exponent; multifractional process; wavelet coefficients; parametric estimation}
\end{abstract}
\section{Introduction}

Multifractional Brownian motion (mBm) is considered as one of the most natural extensions
of fractional Brownian motion (fBm). Nowadays applications of mBm are numerous and growing. Similar to fBm,
mBm has been used in such diverse areas as
geology, image analysis, signal processing, traffic networks and mathematical finance. For instance, we refer to L\'evy-V\'ehel (1995), Bertrand et al. (2012) and Bianchi et al. (2013). In
this brief introduction, we focus on applications to mathematical finance,
which we know best. Since generally neither fBm nor mBm are semi-martingales, Rogers
(1997) pointed out that there would be arbitrage in a market where stocks
are modelled by fBm. However, Cheridito (2003) showed that, if one relaxes
the definition of arbitrage, fBm is an excellent candidate to model
long-term memory in stock markets. Bayraktar et al. (2006)
obtained fBm as the limit of the stock price in an agent-based model where
investors display inertia. Moreover, unlike stock prices, several processes, like
stochastic volatility, exchange rates, or short interest rates do not need
to be semi-martingales for a mathematical model to be arbitrage-free in a
strict sense. For each of these processes, there is empirical evidence of
long-term memory. We refer to Corlay et al. (2014) for stochastic volatility, Xiao et al. (2010) for exchange rates, and Ohashi (2009) for interest rates. Making the
Hurst parameter time-dependent allows to model different regimes of the
stochastic process of interest. For example, in times of financial crisis, asset
volatility rises significantly. Likewise, empirical evidence shows that
there has been periods of different volatility in either exchange rates or
interest rates. This phenomena motivates one to introduce mBm into finance, since unlike fBm, the local regularity of volatilities driven by mBm allows to change via different periods. Let $\{B_{H(t)}(t)\}_t$ denote an mBm with Hurst function $H$. We consider a general model $Y(t)=\Phi(\theta(t)B_{H(t)}(t))$ with $\Phi\in C^2(\mathbb R)$, $\theta\in C^2([0,1])$. In this paper we are interested in estimating $H$, starting from the observations of $Y$. An advantage of our model and methodology is that the
functions $\Phi$ and $\theta$ do not need to be known a priori. This is for instance the case of stochastic
volatility, where the volatility $Y(t)$ is an unknown $C^2$ class function of $X(t)=\theta(t)B_{H(t)}(t)$.

We define the mBm $\{ B_{H(t)}(t)\} _{t\in [0,1]}$ through its harmonizable representation (see Benassi et al. 1997): for $t\in[0,1]$,
\begin{equation}
\label{mBm}
B_{H(t)}(t)=\int_{\mathbb R}\frac{e^{it\xi}-1}{|\xi|^{H(t)+1/2}}{\,\mathrm{d}} \widetilde{W(\xi)},
\end{equation}
where:
\begin{description}
  \item[$-$] The Hurst functional parameter $H(\cdot)$ is  a $\beta$-H\"olderian function with some $\beta\in(0,1)$;
  \item[$-$] The complex-valued stochastic measure ${\,\mathrm{d}}\widetilde{W}$ is defined as the Fourier transform of the real-valued Brownian measure ${\,\mathrm{d}} W$. More precisely, for all $f\in L^2(\mathbb R)$,
      $$
      \int_{\mathbb R}\widehat{f}(t){\,\mathrm{d}}\widetilde{ W(t)}=\int_{\mathbb R}f(t){\,\mathrm{d}} W(t),
      $$
      where $\widehat{f}$ denotes the Fourier transform of $f$:
      $$
      \widehat{f}(\xi)=\int_{\mathbb R}e^{-i\xi t}f(t){\,\mathrm{d}} t.
      $$
\end{description}
Statistically, the most significant feature of mBm, its local H\"older regularity, can be measured by its pointwise H\"older exponent. Recall that for a continuous nowhere differentiable process $\{Y(t)\}_t$, its pointwise H\"older
exponent $\rho_Y$ is a stochastic process defined by: for each $t_0$,
$$
\rho _{Y}(t_0)=\sup\Big\{\alpha\in[0,1]:~\limsup\limits_{\epsilon\rightarrow0}\frac{|Y(t_0+\epsilon)-Y(t_0)|}{|\epsilon|^{\alpha}}=0\Big\}.
$$
Ayache and L\'evy-V\'ehel (2004) show that for each $t_0\in(0,1)$, the pointwise H\"older exponent of $B_{H(t)}(t)$ at each point $t_0$, $\rho _{B_{H}}(t_0)$, is with probability $1$ equal to its Hurst functional
parameter $H(t_0)$: $$\mathbb{P}\big( \rho _{B_{H}}( t_{0}) =H( t_{0})\big )=1.$$

Since generating the scenarios of fBm and mBm relies only on their Hurst
parameters (see Ayache and L\'evy-V\'ehel 2004), the problem of estimation of the Hurst parameter has significant interests. Several results have already been
obtained recently on the estimation of the mBm's Hurst functional
parameter. We refer to Rosenbaum (2008), Coeurjolly (2005, 2006), Bardet
(2002), Bardet and Surgailis (2013) and Bertrand et al. (2013).

Rosenbaum (2008) and Bardet (2002) used wavelet-based methods to study inference problem on fBm, i.e., the Hurst parameter is constant; Coeurjolly (2005, 2006) and Bertrand et al. (2013) studied the estimation of $%
H( t_{0})$  by using respectively generalized variations (see also Chan et al. 1995) and increment ratio statistic method, where a discretized sample path of the mBm is assumed to be observed. Bardet and Surgailis (2013) developed a nonparametric estimation method (based on the increment ratio) for evaluating the local Hurst function of a multifractional Gaussian process (whose increments are asymptotically a multiple of an fBm) which extends mBm, starting from a discrete sample path of this process. In this work we consider a model more general than mBm but slightly less general than the one considered in Bardet and Surgailis (2013). We study a different statistical setting (either the wavelet coefficients of $Y$ or a discrete path of $2^n$ points is observed), by applying a wavelet-based approach. Note that the main advantage of this statistical setting is that, it could be applied to inferential problems when only a set of wavelet coefficients are available (we refer to Delbeke and Van (1995), Abry et al. (2002), Abry and Con\c{c}alv\`es (1997) and the references therein). By a study on the fine regularity property, our central limit theorem provides an explicit form for the limit covariance matrix.  It is worth noting that, the construction of our estimator relies on a sharp estimation of the covariance structure of the mBm's wavelet coefficients. The main technical difficulty arises due to the fact that the covariance structure of mBm is much more complicated than that of fBm. The techniques used to identify such covariance structure is not obvious and thus has its specific interests in statistical inferences on multifractional processes as well as on stochastic analysis.

\subsection{Statistical setting}
Consider the following system: for $t\in[0,1]$,
\begin{equation}
\label{Model1}
Y(t)=\Phi(X(t)),~\mbox{with}~X(t)=\theta(t)B_{H(t)}(t),
\end{equation}
where
\begin{description}
\item[$-$] $\{B_{H(t)}(t)\}_{t\ge0}$ is an mBm defined in (\ref{mBm}). Assume that its Hurst functional parameter $H$ belongs to $C^{2}([0,1])$ and $$\Big[\inf\limits_{t\in[0,1]}H(t),\sup\limits_{t\in[0,1]}H(t)\Big]\subset(0,1);$$
  \item[$-$]  The closed forms of the deterministic functions $\Phi$ and $\theta$ are unknown. However we assume that\begin{description}
   \item[$(1)$] $\Phi\in C^2(\mathbb R)$, $\Phi'\neq 0$ almost everywhere and there exist two constants $0<c_1\le c_2$ such that  $c_1\leq |\Phi'(x)|\leq c_2$ for all $x\in\mathbb R$.
       \item[$(2)$] $\theta\in C^{2}([0,1])$ and $\theta\ne0$ almost everywhere.
       \end{description}
  \item[$-$] Suppose that a discrete trajectory of $Y$: $\big\{Y(u/2^n):u=0,\ldots,2^n\big\}$ is observed for some $n\in\mathbb N$ large enough.
\end{description}
 Our major goal is to evaluate the functional parameter $H(\cdot)$. As in Coeurjolly (2005, 2006), we introduce a pointwise estimation method, namely, the function $H(\cdot)$ is estimated pointwisely for any $t_0\in(0,1)$. Once the time $t_0$ is fixed, the problem becomes a parametric estimation problem. Peng (2011a) studied this problem when $Y(t)$ is some stationary increment process (with $H\equiv\alpha$), where the optimal convergence rate $n^{-1/2}$  is obtained by using the observations $\{Y(k/n):k=0,\ldots,n\}$. However, when $H$ varies via time, the estimation of $H(t_0)$ only relies on the sample size of the observed data in the neighborhood of $Y(t_0)$ and this neighborhood's radius' convergence speed. Hence the convergence speed of the corresponding estimator would be reasonably slower than $n^{-1/2}$ (see e.g. Coeurjolly (2005, 2006)  for a particular case when $Y$ is mBm). In this work, heuristically speaking, since the sample size of the neighborhood data of $Y(t_0)$ for estimating each $H(t_0)$ is about $2^{n+1}\epsilon_n$, it is then believed that a good estimator should have its convergence rate near $(2^{n}\epsilon_n)^{-1/2}+\epsilon_n$. Subject to this statistical setting, we try to get the \enquote{optimal} rate of convergence estimator of $H(t_0)$ by using wavelet basis.
\subsection{Methodology and technical assumptions}
Let the integer $Q\ge1$ and let us pick any mother wavelet $\psi\in C^1([0,1])$  whose first $Q-1$ moments are vanishing:
$$
\int_0^1t^p\psi(t){\,\mathrm{d}} t=0,~\mbox{for}~p=0,\ldots,Q-1;
$$
and
\begin{equation}
\label{wavelet}
\int_0^1t^Q|\psi(t)|{\,\mathrm{d}} t<+\infty~\mbox{and}~\int_0^1t^Q|\psi(t)|{\,\mathrm{d}} t\neq0.
\end{equation}
$Q$ is called the cancellation order of $\psi$. Fix an integer $j\in\{1,2,\ldots,J_n\}$, with $(J_n)_{n\ge1}$ being a subsequence of $\mathbb N$. The wavelet coefficients of $Y(t)$ are given by: for $k\in\{0,\ldots,2^j-1\}$,
\begin{equation}
\label{d}
d_Y(2^{-j},k)=2^{j/2}\int_{0}^{2^{-j}}\psi(2^jt)Y(t+2^{-j}k){\,\mathrm{d}} t.
\end{equation}
We introduce a set of indices corresponding to a neighborhood of each $t_0\in(0,1)$:
$$
\nu_{t_0,2^j}:=\big\{k\in\{0,\ldots,2^j-1\}: |t_0-k2^{-j}|\leq \epsilon_j\big\},
$$
where the radius of this neighborhood $\epsilon_j$ satisfies $\epsilon_j\rightarrow 0$ and $2^j\epsilon_j\rightarrow+\infty$ as $j\rightarrow+\infty$. Then the quadratic form estimator of $Var(Y(t_0))$ is given as:
$$
V_{Y,t_0,j}=\sum_{k\in\nu_{t_0,2^j}}d_Y(2^{-j},k)^2.
$$
For showing the consistency of our estimator and existence of central limit results, we need to  impose  the following technical assumptions to  $(\epsilon_j)_{j\ge1}$, $Q$ and $H(\cdot)$:
\begin{description}
\item[\textbf{(A0):}] $Q\ge2$. This condition leads to $4Q-4\sup_{t\in[0,1]}H(t)>1$. We remark that the latter inequality still holds when $Q=1$ and $\sup_{t\in[0,1]}H(t)<3/4$, as a consequence our results are still valid. However, we won't consider the latter assumptions, because practitioners in data analysis wouldn't expect the unknown functional parameter $H(\cdot)$ should be valued less than $3/4$.
  \item[\textbf{(A1):}] A lower bound of the convergence rate of $\epsilon_j$ toward 0 is given as: $$\lim\limits_{j\rightarrow+\infty}2^{jH(t_0)}\epsilon_j^{2H(t_0)} j^{1/2}|\log\epsilon_j|=0.$$
       \item[\textbf{(A2):}] An upper bound of the convergence rate of $\epsilon_j$ toward 0 is described as: for any $\delta>0$ arbitrarily small, $$\sum\limits_{j=1}^{+\infty}(2^j\epsilon_j)^{-\delta}<+\infty.$$
       \item[\textbf{(A3):}] There exists a constant $c_0>0$ such that $\lim_{j\rightarrow+\infty}\epsilon_{j+1}/\epsilon_j=c_0$.
   \item[\textbf{(A4):}] This condition is stronger than \textbf{(A1)}: $$\lim\limits_{j\rightarrow+\infty}2^{j(H(t_0)+1/2)}\epsilon_j^{2H(t_0)+1/2} j^{1/2}|\log\epsilon_j|=0.$$
\end{description}
Through this paper we assume that \textbf{(A0)} is always satisfied. We will successively show that our estimator  of $H(t_0)$ is weakly consistent under assumption \textbf{(A1)}; it is strongly consistent under assumptions \textbf{(A1)-(A2)}; and it has an asymptotic Gaussian behavior subject to assumptions \textbf{(A2)-(A4)}. Before stating these main results, we make some notation conventions:
\begin{definition}
\begin{description}

\item[(a)] A sequence $(U_n)_{n\in\mathbb N}$ of real-valued random variables is said to be bounded almost surely if and only if
    $$
    \mathbb P\big(\limsup_{n\rightarrow+\infty}|U_n|<+\infty\big)=1;
    $$
    it is said to be bounded in probability if and only if
    $$
    \lim_{\eta\rightarrow+\infty}\limsup_{n\rightarrow+\infty}\mathbb P(|U_n|>\eta)=0.
    $$
    Let $(x_n)_{n\in\mathbb N}$ be a sequence of strictly positive real numbers and let $(U_n)_{n\in\mathbb N}$ be a sequence of positive random variables. The notations
        $$
         U_n=\mathcal O_{a.s.}(x_n)~\mbox{and}~
        U_n=\mathcal O_{\mathbb P}(x_n)
        $$
        mean respectively that the sequence $(x_n^{-1}U_n)_{n\in\mathbb N}$ is bounded almost surely and in probability.
        \item[(b)] We say two sequences of strict positive values $(x_n)_{n\in\mathbb N}$ and $(y_n)_{n\in\mathbb N}$ are equivalent as $n$ tends to infinity if and only if there exist two constants $0<c_1\leq c_2$ such that
            $$
            c_1\leq \liminf_{n\rightarrow+\infty}\frac{x_n}{y_n}\leq \limsup_{n\rightarrow+\infty}\frac{x_n}{y_n}\leq c_2.
            $$
            We denote the relation of equivalence by $x_n\sim y_n$.
            \item[(c)] We use the notations \enquote{$\xrightarrow[]{a.s.}$}, \enquote{$\xrightarrow[]{\mathbb P}$} and \enquote{$\xrightarrow[]{dist}$} to respectively denote convergence $\mathbb P$-almost surely, convergence in probability and convergence in law.
\end{description}
\end{definition}
It is also useful to briefly introduce the steps which lead to the construction of the estimators:
\begin{description}
\item[\textbf{Step 1:}] Identification of $V_{Y,t_0,j}$ starting from the observations $\big\{Y(u/2^n):u=0,\ldots,2^n\big\}$. In order to estimate $V_{Y,t_0,j}$, it suffices to make an identification of $d_Y(2^{-j},k)$ given in (\ref{d}). Such an identification can be naturally obtained by discretization of integrals to the sum as
    \begin{equation}
    \label{defdYn}
\widehat{d_{Y,n}}(2^{-j},k)=2^{j/2}\sum_{l=0}^{2^{n-j}-1}Y(l2^{-n}+k2^{-j})\int_{\frac{l}{2^n}}^{\frac{l+1}{2^n}}\psi(2^{j}t){\,\mathrm{d}} t.
    \end{equation}
    %\Big(\approx2^{j/2}\sum_{l=0}^{2^{n-j}-1}\int_{\frac{l}{2^n}}^{\frac{l+1}{2^n}}\psi(2^jt)Y(t+k2^{-j}){\,\mathrm{d}} t=d_Y(2^{-j},k)\Big).
    By L\'evy's modulus of continuity theorem for mBm (see e.g. Theorem 1.7 in Benassi et al. 1997), $\theta\in C^2([0,1])$, (\ref{supX}) and the fact that
    $$
    |s-t|^{\inf_{u\in[s,t]}H(u)}=\mathcal O\big(|s-t|^{\max\{H(s),H(t)\}}\big)~\mbox{as}~|s-t|\rightarrow0,
    $$
     there exists a positive-valued random variable $C$ all of whose moments are finite, such that for all $s,t\in[0,1]$ with $s\neq t$,
    \begin{eqnarray}
    \label{AH}
 |X(s)-X(t)|&\le&|(\theta(s)-\theta(t))B_{H(s)}(s)|+|\theta(t)(B_{H(s)}(s)-B_{H(t)}(t))|\nonumber\\
 &\le& C|s-t|^{\max\{H(s),H(t)\}} |\log|s-t||^{1/2}.
    \end{eqnarray}
    Then,  the following important relations hold (the proofs are given in the appendix):
    \begin{eqnarray}
    \label{dnds}
    &&|\widehat{d_{Y,n}}(2^{-j},k)-d_Y(2^{-j},k)|=\mathcal O_{a.s.}\big( 2^{-nH(2^{-j}k)-j/2}n^{1/2}\big);\\
    \label{dnd}
    &&\mathbb E|\widehat{d_{Y,n}}(2^{-j},k)-d_Y(2^{-j},k)|^r=\mathcal O\big( 2^{-r(nH(2^{-j}k)+j/2)}\big),~\mbox{for}~r\ge1,
    \end{eqnarray}
    which show that $\widehat{d_{Y,n}}(2^{-j},k)$ is a good estimator of $d_Y(2^{-j},k)$. Next we set
    \begin{equation}
    \label{defdVYn}
    \widehat{V_{n,t_0,j}}=\sum_{k\in\nu_{t_0,2^j}}\widehat{d_{Y,n}}(2^{-j},k)^2,
    \end{equation}
    and show that $\widehat{V_{n,t_0,j}}$ satisfies
    \begin{equation}
    \label{wideV}
    \mathbb E|\widehat{V_{n,t_0,j}}-V_{Y,t_0,j}|^2=\mathcal O(2^{-2(n+j)H(t_0)}\epsilon_j^2).
    \end{equation}
   Please see the appendix for the proof of (\ref{wideV}).
\item[\textbf{Step 2:}] Identification of $H(t_0)$ starting from $V_{X,t_0,j}$. We use similar computations appeared in Peng (2011a). The main result is the following (see Proposition \ref{VarVX} for a more explicit formula):
\begin{equation*}
Var\big(V_{X,t_0,j}\big)=c2^{-j(2H(t_0)+1)}card(\nu_{t_0,2^j})+o(2^{-j(2H(t_0)+1)}card(\nu_{t_0,2^j})),
\end{equation*}
where $c>0$ is a constant which does not depend on $j$; $card(\nu_{t_0,2^j})$ denotes the cardinal of $\nu_{t_0,2^j}$. It is observed that
$$
\lim_{j\rightarrow+\infty}\frac{card(\nu_{t_0,2^j})}{ 2^{j+1}\epsilon_j}=1.
$$
Therefore subject to feasible choices of the sequence $(\epsilon_j)_j$, the following convergence can accordingly take place in probability or almost surely:
\begin{equation}
\label{as1'}
\lim_{j\rightarrow+\infty}\frac{V_{X,t_0,j}}{c2^{-2jH(t_0)}\epsilon_j}=1,
\end{equation}
where $c>0$ is some constant not depending on $j$.
\item[\textbf{Step 3:}] Identification of $H(t_0)$ starting from $\{\widehat{V_{n,t_0,J_n}}\}_{n\ge1}$.

Under assumption \textbf{(A1)} (resp. \textbf{(A1)-(A2)}), one has the following relation of equivalence (see (\ref{as1'}), (\ref{Taylor2}),  (\ref{DiffV_p}) and (\ref{DiffV_as})): for $J_n\in\{0,\ldots,n-1\}$,
$$
\widehat{V_{n,t_0,J_n}}\sim V_{Y,t_0,J_n}\sim V_{X,t_0,J_n}\sim 2^{-2J_nH(t_0)}\epsilon_{J_n},~\mbox{as $J_n\rightarrow+\infty$}
$$
in probability (resp. almost surely). As a consequence,
$$
\widehat{H_{Y,2^{J_n},n}}(t_0)=\frac{\log\left(\frac{\epsilon_{J_n+1}}{\epsilon_{J_n}}\right)-\log(\frac{\widehat{V_{n,t_0,J_n+1}}}{\widehat{V_{n,t_0,J_n}}})}{2\log2}
$$
is a consistent  estimator of $H(t_0)$. We remark that the speed of convergence relies on the choice of $J_n$. More details on the choice of $J_n$ will be discussed in Theorem \ref{VYdn}.
\end{description}
\section{Some preliminary results}
\subsection{Identification of the covariance structure of the wavelet coefficients of $X$}
In this part we provide a sharp estimation of the covariance structure of wavelet coefficients of $X(t)=\theta(t)B_{H(t)}(t)$. For $a>0$, define the wavelet coefficient $d_X\left( a,k\right) $ of $X$ by
\begin{equation}
\label{da}
d_X\left( a,k\right) =\frac{1}{\sqrt{a}}\int_{0}^{a}\psi \big( \frac{t}{a}%
\big) X(t+ak){\,\mathrm{d}} t,
\end{equation}
for all $k\in \left\{0,1,...,\left[ a^{-1}\right]-1 \right\},$ with $[\cdot]$ being the integer part function.\\
The following proposition provides a fine identification of  $d_X\left( a,k\right)$'s covariance structure .
\begin{proposition}
\label{prop1}
Let $a,b>0$ satisfy that, there exist two constants $c_1,c_2>0$ such that $0<c_1\leq \liminf\limits_{a,b\rightarrow0}\frac{a}{b}\leq\limsup\limits_{a,b\rightarrow0}\frac{a}{b}\leq c_2$\footnote{We will denote this relation by $a\sim b$.}.  Let $k\in\{0,\ldots,[a^{-1}]-1\}$, $k'\in\{0,\ldots,[b^{-1}]-1\}$.\\
$(1)$  If $ak\neq bk'$ and $\sup\limits_{s,t\in[0,1]}|\frac{at-bs}{ak-bk'}|\leq1$,  then we have
\begin{eqnarray}
\label{prop11}
&&Cov\left( d_X\left( a,k\right) ,d_X\left( b,k^{\prime }\right) \right) =\frac{%
C_1(H(ak)+H(bk'),Q,a/b)\theta(ak)\theta(bk')(ab) ^{Q+1/2}}{\left\vert ak-bk^{\prime }\right\vert
^{2Q-H\left( ak\right) -H\left( bk^{\prime }\right) }}\nonumber\\
&&+(ab)^{(H(ak)+H(bk')+1)/2}T(k,k',Q,a,b),
\end{eqnarray}%
where the closed form of $C_1(\cdot,\cdot,\cdot)$ is given in $(\ref{computeI00})$ and the remaining terms $\{T(k,k',Q,a,b)\}_{k,k'}$ verify
$$
\sum_{0\leq k\leq [a^{-1}]-1,0\leq k'\leq [b^{-1}]-1,ak\neq bk'}\!\!\!\!\!\!\!\!\!\!\!\!\!\!\big(T(k,k',Q,a,b)\big)^2=\mathcal O((ab)^{1/2}\log a\log b),~\mbox{as $a,b\rightarrow0$}.
$$
$(2)$ If  $b/a=k/k'=\varrho>0$, we have
\begin{equation}
\label{prop12}
Cov\big( d_X\left( a,k\right),d_X(b,k')\big)=C_2(ak,\varrho)a^{2H\left( ak\right) +1}+\mathcal{O%
}\big( a^{2H\left( ak\right) +2}|\log a|^4\big),
\end{equation}
where the explicit form of $C_2(\cdot,\cdot)$ is given in $(\ref{c2})$.
\end{proposition}
The proof of Proposition \ref{prop1} is given in the appendix. The following lemma is a straightforward consequence of Proposition \ref{prop1} (it suffices to take $a=b=2^{-j}$ in Proposition \ref{prop1}.) that we will rely on heavily through the remaining context.
\begin{lemma}
\label{lemma1}
For $k,k'\in\{0,\ldots,2^{j}-1\}$ and $k\neq k'$, we have
\begin{eqnarray}
\label{lemma11}
&&Cov\left( d_X( 2^{-j},k) ,d_X( 2^{-j},k^{\prime }) \right)=\frac{%
c_1(2^{-j}k,2^{-j}k')2 ^{-j(1+H(2^{-j}k)+H(2^{-j}k'))}}{|k-k'|^{2Q-H(2^{-j}k) -H(2^{-j}k') }}\nonumber\\
&&+2 ^{-j(1+H(2^{-j}k)+H(2^{-j}k'))}T(k,k',Q,2^{-j},2^{-j}),
\end{eqnarray}
where $c_1(\cdot,\cdot)\in C^2([0,1])$ is defined as
\begin{equation}
\label{c1}
c_1(x,y)=C_1(H(x)+H(y),Q,1)\theta(x)\theta(y)~\mbox{in (\ref{computeI00})}.
\end{equation}
And
\begin{equation}
\label{lemme12}
\mathbb E| d_X( 2^{-j},k)|^{2}=c_2(2^{-j}k)2^{-j(2H(2^{-j}k) +1)}+\mathcal{O%
}( 2^{-2j(H(2^{-j}k) +1)}j^4),
\end{equation}
where $c_2(\cdot)\in C^2([0,1])$ is defined as $c_2(x)=C_2(x,1)$ in (\ref{c2}).
\end{lemma}
\subsection{Identification of $H(t_0)$ when $\{d_X(2^{-j},k)\}_{k=0,\ldots,2^{j}-1}$ are observed}
In this section we first construct a consistent estimator of $H(t_0)$, when the wavelet coefficients of $X$: $\{d_X(2^{-j},k)\}_{k=0,\ldots,2^j-1}$ can be straightforwardly observed (i.e. $\Phi\equiv Id$).  To this end we let
\begin{equation}
\label{VXV}
V_{X,t_0,j}=\sum_{k\in\nu_{t_0,2^j}}d_X(2^{-j},k)^2.
\end{equation}
The following proposition provides sharp identifications of $V_{X,t_0,j}$'s first order and second order moments.
\begin{proposition}
\label{VarVX}
For any $t_0\in(0,1)$,
\begin{equation}
\label{VarV'}
\mathbb E\big(V_{X,t_0,j}\big)=2c_2(t_0)2^{-2jH(t_0)}\epsilon_j+\mathcal O\big(2^{-2jH(t_0)}(2^{-j}+j\epsilon_j^2+2^{-j}j^4\epsilon_j)\big);
\end{equation}
and
\begin{eqnarray}
\label{VarV}
&&Var\big(V_{X,t_0,j}\big)=c_3(t_0)2^{-j(4H(t_0)+1)}\epsilon_j+\mathcal O(2^{-j(4H(t_0)+1)}(2^{-j}\!\!+j\epsilon_j^2+2^{-j/2}j\epsilon_j)\big),\nonumber\\
\end{eqnarray}
where the constants $c_1(t_0,t_0)$ and $c_2(t_0)$ are given in (\ref{computeI00}) and (\ref{c2}) respectively, and
\begin{equation}
\label{c3t0}
c_3(t_0)=4\big(c_2(t_0)^2+c_1(t_0,t_0)^2\sum\limits_{l\in\mathbb Z,l\neq0}|l|^{4H(t_0)-4Q}\big).
\end{equation}
\end{proposition}
The proof of Proposition \ref{VarVX} is given in the appendix.

We remark that, since a sequence of random variables bounded almost surely is also bounded in probability, therefore the following identification of $V_{X,t_0,j}$ holds, thanks to Proposition \ref{VarVX}, Chebyshev's inequality and Borel-Cantelli's lemma.
\begin{lemma}
\label{IV}
If $\epsilon_j=\mathcal O(j^{-1})$, we have
\begin{eqnarray}
\label{IV1}
V_{X,t_0,j}&=&2c_2(t_0)2^{-2jH(t_0)}\epsilon_j+\mathcal O_{a.s.}\big(2^{-2jH(t_0)}(j\epsilon_j^2+2^{-j}j^4\epsilon_j)\big)\nonumber\\
&&+\mathcal O_{\mathbb P}(2^{-j(2H(t_0)+1/2)}\epsilon_j^{1/2}).
\end{eqnarray}
As a consequence,
\begin{description}
\item[(a)] if $\lim\limits_{j\rightarrow+\infty}j\epsilon_j=0$, then
\begin{equation}
\label{IV2}
\frac{V_{X,t_0,j}}{2c_2(t_0)2^{-2jH(t_0)}\epsilon_j}-1=\mathcal O_{\mathbb P}(j\epsilon_j+2^{-j}j^4+(2^j\epsilon_j)^{-1/2});
\end{equation}
\item[(b)] if $\lim\limits_{j\rightarrow+\infty}j\epsilon_j=0$ and $\sum\limits_{j=1}^{+\infty}(2^j\epsilon_j)^{-\delta}<+\infty$ for arbitrarily small $\delta>0$, then
\begin{equation}
\label{IV3}
\frac{V_{X,t_0,j}}{2c_2(t_0)2^{-2jH(t_0)}\epsilon_j}-1=\mathcal O_{a.s.}(j\epsilon_j+2^{-j}j^4+(2^j\epsilon_j)^{-1/2+\delta}),
\end{equation}
for arbitrarily small $\delta>0$.
\end{description}
\end{lemma}
\textbf{Proof.} By using Chebyshev's inequality, (\ref{VarV})  and the condition that $\epsilon_j=\mathcal O(j^{-1})$, we get for any $\eta>0$,
\begin{eqnarray}
\label{IV4}
&&\mathbb P\Big(2^{j(2H(t_0)+1/2)}\epsilon_j^{-1/2}\big|V_{X,t_0,j}-\mathbb E(V_{X,t_0,j})\big|\ge\eta\Big)\leq\frac{ 2^{j(4H(t_0)+1)}\epsilon_j^{-1} Var (V_{X,t_0,j})}{\eta^2}\nonumber\\
&&\leq c\frac{2^{j(4H(t_0)+1)}2^{-j(4H(t_0)+1)}}{\eta^2}= \frac{c}{\eta^2},
\end{eqnarray}
where $c>0$ is some constant which does not depend on $j$. This implies
\begin{equation}
\label{IV5}
V_{X,t_0,j}=\mathbb E(V_{X,t_0,j})+\mathcal O_{\mathbb P}(2^{-j(2H(t_0)+1/2)}\epsilon_j^{1/2}).
\end{equation}
Then it follows from (\ref{VarV'}), (\ref{IV5}) and the fact that $\lim_{j\rightarrow+\infty}2^{-j/2}\epsilon_{j}^{-1/2}=0 $ that
\begin{eqnarray*}
&&V_{X,t_0,j}=2c_2(t_0)2^{-2jH(t_0)}\epsilon_j+\mathcal O_{a.s.}\big(2^{-2jH(t_0)}(2^{-j}+j\epsilon_j^2+2^{-j}j^4\epsilon_j)\big)\nonumber\\
&&~~+\mathcal O_{\mathbb P}(2^{-j(2H(t_0)+1/2)}\epsilon_j^{1/2})\nonumber\\
&&=2c_2(t_0)2^{-2jH(t_0)}\epsilon_j+\mathcal O_{a.s.}\big(2^{-2jH(t_0)}(j\epsilon_j^2+2^{-j}j^4\epsilon_j)\big)\nonumber\\
&&~~+\mathcal O_{\mathbb P}(2^{-j(2H(t_0)+1/2)}\epsilon_j^{1/2}).
\end{eqnarray*}
(\ref{IV1}) has been proven. Note that (\ref{IV2}) follows straightforwardly from (\ref{IV1}). Now we only need to show (\ref{IV3}) holds. From (\ref{IV5}) and Chebyshev's inequality, we observe that there exists a constant $c>0$ which does not depend on $\eta$ nor on $j$ such that for any $\eta>0$,
\begin{equation}
\label{IV6}
\mathbb P\Big((2^j\epsilon_j)^{(1-\delta)/2}\big|\frac{V_{X,t_0,j}}{\mathbb E(V_{X,t_0,j})}-1\big|>\eta\Big)\leq \frac{c(2^j\epsilon_j)^{-\delta}}{\eta^2}.
\end{equation}
Since $\sum\limits_{j=1}^{+\infty}(2^j\epsilon_j)^{-\delta}<+\infty$,  then applying Borel-Cantelli's lemma leads to
$$
\frac{V_{X,t_0,j}}{\mathbb E(V_{X,t_0,j})}-1=\mathcal O_{a.s.}((2^j\epsilon_j)^{-1/2+\delta/2}).
$$
Further observe that
$$
\mathbb E(V_{X,t_0,j})=\mathcal O_{a.s.}(2^{-2jH(t_0)}\epsilon_j).
$$
Therefore (\ref{IV3}) follows. Lemma \ref{IV} has been proven. $\square$

Now we are ready to state the main results of this section. The following theorem constructs a consistent estimator of $H(t_0)$ starting from the wavelet coefficients $\{d_X(2^{-j},k)\}_{k=0,\ldots,2^{j}-1}$.
\begin{theorem}
\label{VXd}
For $t_0\in(0,1)$, denote by
$$
\widehat{H_{X,2^j}}(t_0)=\frac{\log\big(\frac{\epsilon_{j+1}}{\epsilon_j}\big)+\log\big(\frac{V_{X,t_0,j}}{V_{X,t_0,j+1}}\big)}{2\log 2}.
$$
\begin{description}
\item[(a)] If $\lim\limits_{j\rightarrow+\infty}j\epsilon_j=0$, then
\begin{equation}
\label{VXd2}
\widehat{H_{X,2^j}}(t_0)-H(t_0)=\mathcal O_{\mathbb P}(j\epsilon_j+2^{-j}j^4+(2^j\epsilon_j)^{-1/2});
\end{equation}
\item[(b)] if $\lim\limits_{j\rightarrow+\infty}j\epsilon_j=0$ and $\sum\limits_{j=1}^{+\infty}(2^j\epsilon_j)^{-\delta}<+\infty$ for any $\delta>0$ arbitrarily small, then
    \begin{equation}
\label{VXd3}
\widehat{H_{X,2^j}}(t_0)-H(t_0)=\mathcal O_{a.s.}(j\epsilon_j+2^{-j}j^4+(2^j\epsilon_j)^{-1/2+\delta}).
\end{equation}
\end{description}

\end{theorem}
\textbf{Proof.} This theorem is a straightforward consequence of Lemma \ref{IV}, as we observe that, by construction, $\widehat{H_{X,2^j}}(t_0)$ verifies
\begin{eqnarray*}
&&\widehat{H_{X,2^j}}(t_0)-H(t_0)\\
&&=\frac{1}{2\log 2}\Big(\log\big(\frac{V_{X,t_0,j}}{2c_2(t_0)2^{-2jH(t_0)}\epsilon_{j}}\big)-\log\big(\frac{V_{X,t_0,j+1}}{2c_2(t_0)2^{-2(j+1)H(t_0)}\epsilon_{j+1}}
\big)\Big).
\end{eqnarray*}
Then, by Lemma \ref{IV}, the fact that the logarithm function belongs to $C^1((0,+\infty))$, continuous mapping theorem and the mean value theorem, we get Theorem \ref{VXd}.  $\square$

In Theorem \ref{VXd}, if we assume $\epsilon_j=j^{\alpha}2^{-j\gamma}$ with some $\alpha\in\mathbb R$ and $\gamma\in(0,1)$, then by elementary computation we obtain the following corollary, which  leads to choose $\epsilon_j\approx j^{-2/3}2^{-j/3}$ to obtain the best rate of convergence of the estimator $\widehat{H_{X,2^j}}(t_0)$. This result is similar to the choice of $b(N)$ in Bardet (2002) in the setting of fBm.
\begin{corollary}
\label{cor}
Under assumption $\epsilon_j=j^{\alpha}2^{-j\gamma}$ with $\alpha\in\mathbb R$ and $\gamma\in(0,1)$,
\begin{description}
\item[(a)] Taking $\alpha=-2/3$ and $\gamma=1/3$, we have
$$
\sup_{t\in(0,1)}|\widehat{H_{X,2^j}}(t)-H(t)|=\mathcal O_{\mathbb P}((j2^{-j})^{1/3}).
$$
\item[(b)] For small $\delta>0$, taking $\alpha=-\frac{2}{3-2\delta}$ and $\gamma=\frac{1-2\delta}{3-2\delta}$, we obtain
$$
\sup_{t\in(0,1)}|\widehat{H_{X,2^j}}(t)-H(t)|=\mathcal O_{a.s.}((j2^{-j})^{\frac{1-2\delta}{3-2\delta}}).
$$
\end{description}
\end{corollary}

The second main result of this part describes  an asymptotic Gaussian behavior of the estimator of $H(t_0)$, where the limit covariance matrix is precisely given, which depends on $H(t_0)$.
\begin{theorem}
\label{cltX}
For $t_0\in(0,1)$, if $\lim\limits_{j\rightarrow+\infty}j\epsilon_j=0$ and assumptions $\textbf{(A2)-(A3)}$ are verified, then
    \begin{equation}
\label{cltX1}
\sqrt{2^{j+1}\epsilon_j}\Big(\widehat{H_{X,2^j}}(t_0)-H(t_0)\Big)\xrightarrow[j\rightarrow+\infty]{dist}\mathcal N\big(0,\tilde{c}(t_0)\big),
\end{equation}
where the constant $$\tilde{c}(t_0)=\frac{1}{(2\log 2)^2}\Big(\big((2c_0)^{-1}+1\big)\frac{c_3(t_0)}{c_2(t_0)^2}-2(2c_0)^{-1/2}c_4(t_0)\Big),$$ with $c_0$ given in $\textbf{(A3)}$ and $c_3(t_0)$ given in (\ref{c3t0}).
\end{theorem}In order to prove Theorem \ref{cltX}, we rely heavily on the following proposition.
\begin{proposition}
\label{prop:clt1}
For any $t_0\in(0,1)$ and any integer $j\ge1$, denote by
$$
U_{t_0,j}=\sqrt{card(\nu_{t_0,2^j})}\Big(\frac{1}{card(\nu_{t_0,2^j})}\sum_{k\in\nu_{t_0,2^j}}\frac{d_X(2^{-j},k)^2}{\mathbb E(d_X(2^{-j},k)^2)}-1\Big).
$$
Then
$$
\big(U_{t_0,j},U_{t_0,j+1}\big)\xrightarrow[j\rightarrow+\infty]{dist}\mathcal N(0,\Sigma),
$$
where $\Sigma=(\sigma_{ij})_{i,j\in\{1,2\}}$ with $\sigma_{11}=\sigma_{22}=c_3(t_0)/(c_2(t_0))^2$ and
\begin{eqnarray}
\label{c4}
&&\sigma_{12}=\sigma_{21}=c_4(t_0):=\frac{(2c_0)^{1/2}}{c_2(t_0)^2}\Big(C_2(t_0,1/2)^22^{2H(t_0)+1}+2c_1(t_0,t_0)^2\nonumber\\
&&~~+2^{2Q-2H(t_0)+1}C_1(H(t_0),H(t_0),Q,2)^2\theta(t_0)^4\!\!\!\sum_{l\in\mathbb Z,|l|\ge2}\!\!\frac{1}{|l|^{4Q-4H(t_0)}}\Big).\nonumber\\
\end{eqnarray}
\end{proposition}
\textbf{Proof.}  Following the similar method as in Bardet (2000), we show that for the empirical average of $d_X( 2^{-j},k)^2/\mathbb E(d_X(2^{-j},k)^2)$, a multivariate central limit theorem holds thanks to  a Lindeberg's condition. More precisely, for $j$ big enough and for $k\in\nu_{t_0,2^j}$, $k'\in\nu_{t_0,2^{j+1}}$, denote by
$$
T_{j,k,k'}=Cov\left( \frac{d_X( 2^{-j},k)^2}{\mathbb E(d_X(2^{-j},k)^2)} ,\frac{d_X( 2^{-{(j+1)}},k')^2}{\mathbb E(d_X(2^{-{(j+1)}},k')^2)} \right),
$$
a Lindeberg's condition thus can be deduced from the following relations:
 \begin{itemize}
  \item For $|2k-k'|=0$,
  \begin{equation}
 \label{clt22'}
 T_{j,k,k'}=2^{2H(2^{-j}k)+2}\Big(\frac{C_2(2^{-j}k,1/2)}{c_2(2^{-j}k)}\Big)^2+\mathcal O(2^{-j}j^4).
 \end{equation}
  \item For $|2k-k'|=1$,
   \begin{equation}
 \label{clt22''}
T_{j,k,k'}=2\Big(\frac{c_1(2^{-j}k,2^{-(j+1)}k')^2}{c_2(2^{-j}k)c_2(2^{-(j+1)}k')}\Big)+\mathcal O\big(T(2k,k',Q,2^{-j},2^{-j})\big).
 \end{equation}
  \item For $|2k-k'|\ge2$,
  \begin{eqnarray}
 \label{clt2}
 T_{j,k,k'}&=&2\Big(\frac{(C_1(H(2^{-j}k)+H(2^{-(j+1)}k'),Q,2)\theta(2^{-j}k)\theta(2^{-(j+1)}k'))^2}{c_2(2^{-j}k)c_2(2^{-(j+1)}k')|2k-k'|^{4Q-2H(2^{-j}k)-2H(2^{-(j+1)}k')}}\Big)\nonumber\\
 &&+\mathcal O\big(T(k,k',Q,2^{-j},2^{-(j+1)})\big).
 \end{eqnarray}
\end{itemize}
 To show (\ref{clt22'})-(\ref{clt2}) hold, we first observe, from (\ref{identcovGaus}), that
 \begin{equation}
 \label{idT}
 T_{j,k,k'}=2\bigg(Cov\Big( \frac{d_X( 2^{-j},k)}{\sqrt{\mathbb E(d_X(2^{-j},k)^2)}} ,\frac{d_X( 2^{-{(j+1)}},k')}{\sqrt{\mathbb E(d_X(2^{-{(j+1)}},k')^2)}} \Big)\bigg)^2.
 \end{equation}
 Therefore (\ref{clt22'}) results from (\ref{idT}), (\ref{prop12}) and (\ref{lemme12}). In order to obtain (\ref{clt22''}), it suffices to take $l=2k$ for $k\in\nu_{t_0,2^j}$, then $l,k'$ belong to $\nu_{t_0,2^{j+1}}$ and they satisfy $\sup_{s,t\in[0,1]}|(t-s)/(l-k')|=\sup_{s,t\in[0,1]}|s-t|\leq1$. This entails that (\ref{prop11}) can be applied on $(a,b,l,k')$ by setting $a=b=2^{-(j+1)}$ and $|l-k'|=1$. As a consequence (\ref{clt22''}) follows from (\ref{idT}), (\ref{lemma11}) and (\ref{lemme12}). For proving (\ref{clt2}), we
  just plug $a=2^{-j}$, $b=2^{-(j+1)}$ into (\ref{prop11}) and then use (\ref{idT}) and (\ref{lemme12}). Using this identification of $T_{j,k,k'}$'s, we show the Lindeberg's condition (the same as in Bardet 2000) is verified and the central limit theorem holds. Now it remains to show
   \begin{eqnarray}
   \label{convcov}
   &&Cov\big(U_{t_0,j},U_{t_0,j+1}\big)\xrightarrow[j\rightarrow+\infty]{}c_4(t_0)~~(\mbox{given in (\ref{c4})});\\
   \label{convvar}
   &&Var(U_{t_0,j})\xrightarrow[j\rightarrow+\infty]{}\frac{c_3(t_0)}{c_2(t_0)^2}.
   \end{eqnarray}
   We only prove (\ref{convcov}) holds since (\ref{convvar}) can be followed by quite a similar way. Remark that $\lim_{j\rightarrow+\infty}\sup_{k\in\nu_{t_0,2^j}}2^{-j}k=t_0$; and for $j$ large enough,
    $$card\big(\{k\in\nu_{t_0,2^j},k'\in\nu_{t_0,2^{j+1}},|2k-k'|=0\}\big)=card(\nu_{t_0,2^j});$$
    $$card\big(\{k\in\nu_{t_0,2^j},k'\in\nu_{t_0,2^{j+1}},|2k-k'|=l\}\big)=2card(\nu_{t_0,2^j})+\mathcal O(1)~\mbox{for}~l\ge1;$$
     and for a function $f$ continuous on $t_0$, $$\lim\limits_{j\rightarrow+\infty}\frac{1}{card(\nu_{t_0,2^j})}\sum\limits_{k\in\nu_{t_0,2^j}}f(2^{-j}k)=f(t_0).$$ Considering all the above facts and (\ref{clt22'})-(\ref{clt2}), (\ref{VarV21})-(\ref{suml}), $card(\nu_{t_0,2^j})\sim  2^{j}\epsilon_j$ and assumption \textbf{(A3)}, we obtain
\begin{eqnarray*}
\label{clt3}
&&Cov\big(U_{t_0,j},U_{t_0,j+1}\big)=\frac{1}{\sqrt{card(\nu_{t_0,2^j})card(\nu_{t_0,2^{j+1}})}}\nonumber\\
&&~~\times\Big(\sum_{(k\in\nu_{t_0,2^j},k'\in\nu_{t_0,2^{j+1}},2k=k')}T_{j,k,k'}
+\sum_{|2k-k'|=1}T_{j,k,k'}\\
&&~~+\sum_{l=2}^{+\infty}\sum_{|2k-k'|=l}T_{j,k,k'}\Big)\xrightarrow[j\rightarrow+\infty]{}c_4(t_0)~~~~\mbox{(given in (\ref{c4}))}.
\end{eqnarray*}
Finally, we have proved Proposition \ref{prop:clt1}. $\square$

Next we present the following classical result (see e.g. Oehlert 1992).
\begin{proposition}[Multivariate delta rule]
\label{deltamethod}
Let the estimators $\{(X_n,Y_n)\}_{n\in\mathbb N}$ (valued in $(0,+\infty)^2$) of $(\theta_1,\theta_2)$ satisfy the following central limit theorem:
$$
h(n)\Big((X_n,Y_n)-(\theta_1,\theta_2)\Big)\xrightarrow[n\rightarrow+\infty]{dist}\mathcal N(0,\Sigma),
$$
where $\Sigma$ denotes the covariance matrix of the limit distribution and $(h(n))_n$ is a sequence of positive numbers tending to infinity. Let $g:~(0,+\infty)^2\rightarrow \mathbb R^p$ ($p=1$ or $2$) belong to $C^1((0,+\infty)^2)$, then the following convergence in law holds:
$$
h(n)\Big(g(X_n,Y_n)-g(\theta_1,\theta_2)\Big)\xrightarrow[n\rightarrow+\infty]{dist}\mathcal N\Big(0,\nabla g(\theta_1,\theta_2)^T\Sigma\nabla g(\theta_1,\theta_2)\Big),
$$
where $\nabla g(\theta_1,\theta_2)^T$ denotes the transpose of the gradient of $g$ on $(\theta_1,\theta_2)$.
\end{proposition}
Note that if $p=2$ in Proposition \ref{deltamethod}, the gradient of $g$ becomes a Jacobian matrix.\\
\textbf{Proof of Theorem \ref{cltX}.} For simplifying notation we denote by
$$
\widehat{U_{t_0,j}}=\frac{1}{2^{j+1}\epsilon_j}\frac{\sum_{k\in\nu_{t_0,2^j}}d_X(2^{-j},k)^2}{c_2(t_0)2^{-j(2H(t_0)+1)}}-1.
$$Therefore the following decomposition holds:
\begin{eqnarray}
\label{dec1}
&&\sqrt{2^{j+1}\epsilon_j}\widehat{U_{t_0,j}}-U_{t_0,j}=\frac{1}{\sqrt{2^{j+1}\epsilon_j}}\sum_{k\in\nu_{t_0,2^j}}d_X(2^{-j},k)^2\nonumber\\
&&~~~~\times\Big(\frac{1}{c_2(t_0)
2^{-j(2H(t_0)+1)}}-\frac{1}{\mathbb E(d_X(2^{-j},k)^2)}\sqrt{\frac{2^{j+1}\epsilon_j}{card(\nu_{t_0,2^j})}}\Big)\nonumber\\
&&~~~~+\big(\sqrt{card(\nu_{t_0,2^j})}-\sqrt{2^{j+1}\epsilon_j}\big)\nonumber\\
&&=\frac{1}{\sqrt{2^{j+1}\epsilon_j}}\sum_{k\in\nu_{t_0,2^j}}d_X(2^{-j},k)^2\nonumber\\
&&~~~~\times\Big(\frac{\mathbb E(d_X(2^{-j},k)^2)\sqrt{card(\nu_{t_0,2^j})}-c_2(t_0)
2^{-j(2H(t_0)+1)}\sqrt{2^{j+1}\epsilon_j}}{c_2(t_0)
2^{-j(2H(t_0)+1)}\mathbb E(d_X(2^{-j},k)^2)\sqrt{card(\nu_{t_0,2^j})}}\Big)\nonumber\\
&&~~~~+\big(\sqrt{card(\nu_{t_0,2^j})}-\sqrt{2^{j+1}\epsilon_j}\big).
\end{eqnarray}
Since the fact that $card(\nu_{t_0,2^j})=2^{j+1}\epsilon_j+\mathcal O(1)$ implies
\begin{equation}
\label{cardIdent}
\sqrt{card(\nu_{t_0,2^j})}=\sqrt{2^{j+1}\epsilon_j}+\mathcal O(2^{-j/2}\epsilon_j^{-1/2}).
\end{equation}
And also by taking $r=4$ in (\ref{VnV2}), we have
\begin{equation}
\label{das}
\mathbb E|d_X(2^{-j},k)|^4=\mathcal O(2^{-j(4H(t_0)+2)}).
\end{equation}
By (\ref{lemme12}) (\ref{dec1}), (\ref{Jensen}), (\ref{cardIdent}) and (\ref{das}), we then obtain
\begin{eqnarray*}
&&\mathbb E\big|\sqrt{2^{j+1}\epsilon_j}\widehat{U_{t_0,j}}-U_{t_0,j}\big|^2=\frac{card(\nu_{t_0,2^j})^2}{card(\nu_{t_0,2^j})}\mathcal O(2^{-j(4H(t_0)+2)})\nonumber\\
\!\!\!\!&&\times\Big(\frac{(c_2(t_0)2^{-j(2H(t_0)+1)}+\mathcal O(2^{-j(2H(t_0)+2)}j^4))(\sqrt{2^{j+1}\epsilon_j}+\mathcal O(2^{-(j+1)/2}\epsilon_j^{-1/2}))}{(c_2(t_0)
2^{-j(2H(t_0)+1)})^2\sqrt{2^{j+1}\epsilon_j}}\nonumber\\
&&-\frac{1}{c_2(t_0)
2^{-j(2H(t_0)+1)}}\Big)^2+\Big(\sqrt{2^{j+1}\epsilon_j}+\mathcal O(2^{-(j+1)/2}\epsilon_j^{-1/2})-\sqrt{2^{j+1}\epsilon_j}\Big)^2\nonumber\\
&&=\mathcal O\big(2^{-j}\epsilon_jj^8+(2^{j}\epsilon_j)^{-1}\big).
\end{eqnarray*}
It follows from Markov's inequality that there exists $c>0$ such that
\begin{eqnarray*}
&&\mathbb P\Big(\big|\sqrt{2^{j+1}\epsilon_j}\widehat{U_{t_0,j}}-U_{t_0,j}\big|>\eta\Big)\le \eta^{-2}\mathbb E\big|\sqrt{2^{j+1}\epsilon_j}\widehat{U_{t_0,j}}-U_{t_0,j}\big|^2\\
&&\le c(2^{-j}\epsilon_jj^8+(2^{j}\epsilon_j)^{-1}).
\end{eqnarray*}
The assumption \textbf{(A2)} then allows us to apply Borel-Cantelli's lemma to obtain
 $$
 \big|\sqrt{2^{j+1}\epsilon_j}\widehat{U_{t_0,j}}-U_{t_0,j}\big|\xrightarrow[]{a.s.}0.
 $$
 Therefore, it follows from Proposition \ref{prop:clt1}  and continuous mapping theorem that
$$
\sqrt{2^{j+1}\epsilon_j}\Big(\widehat{U_{t_0,j}},\sqrt{2\epsilon_{j+1}/\epsilon_j}\widehat{U_{t_0,j+1}}\Big)\xrightarrow[j\rightarrow+\infty]{dist}\mathcal N(0,\Sigma).
$$
Since $\lim_{j\rightarrow+\infty}\sqrt{2\epsilon_{j+1}/{\epsilon_j}}=\sqrt{2c_0}$, using  Slutsky's theorem, we get
$$
\sqrt{2^{j+1}\epsilon_j}\Big(\widehat{U_{t_0,j}},\sqrt{2c_0}\widehat{U_{t_0,j+1}}\Big)\xrightarrow[j\rightarrow+\infty]{dist}\mathcal N(0,\Sigma).
$$
This is in fact equivalent to
\begin{equation}
\label{dec2}
\sqrt{2^{j+1}\epsilon_j}(\widehat{U_{t_0,j}},\widehat{U_{t_0,j+1}})\xrightarrow[j\rightarrow+\infty]{dist}\mathcal N(0,\widetilde{\Sigma}),
\end{equation}
with $\widetilde{\Sigma}=\Big({}_{(2c_0)^{-1/2}\sigma_{21}}^{\sigma_{11}}~{}_{(2c_0)^{-1}\sigma_{22}}^{(2c_0)^{-1/2}\sigma_{12}}\Big)$.
Then by applying Proposition \ref{deltamethod} to (\ref{dec2}) with $g_1(x,y)=(\log(x),\log(y))$,  $(X_j,Y_j)=\big(\widehat{U_{t_0,j}}+1,\widehat{U_{t_0,j+1}}+1\big)$, $(\theta_1,\theta_2)=(1,1)$ and $h(j)=\sqrt{2^{j+1}\epsilon_j}$, we get
\begin{eqnarray}
\label{dec3}
&&\sqrt{2^{j+1}\epsilon_j}\Big(\log(\sum_{k\in\nu_{t_0,2^i}}d_X(2^{-i},k)^2)+2iH(t_0)\log(2)-\log(2c_2(t_0)\epsilon_i)\Big)_{i=j,j+1}\nonumber\\
&&~~\xrightarrow[j\rightarrow+\infty]{dist}\mathcal N(0,\widetilde{\Sigma}),
\end{eqnarray}
because the Jacobian matrix $\nabla g_1(1,1)=Id$.

Now we apply again Proposition \ref{deltamethod} to (\ref{dec3}) with $g_2(x,y)=\frac{x-y}{2\log 2}$, to get
$$
\sqrt{2^{j+1}\epsilon_j}\Big(\widehat{H_{X,2^j}}(t_0)-H(t_0)\Big)\xrightarrow[j\rightarrow+\infty]{dist}\mathcal N\big(0,\tilde{c}(t_0)\big),
$$
where
\begin{equation}
\label{tildec}
\tilde{c}(t_0)=\frac{\big({}_{-1}^{~1}\big)^T\widetilde{\Sigma}\big({}_{-1}^{~1}\big)}{(2\log 2)^2}=\frac{1}{(2\log 2)^2}\Big(\big((2c_0)^{-1}+1\big)\frac{c_3(t_0)}{c_2(t_0)^2}-2(2c_0)^{-1/2}c_4(t_0)\Big),
\end{equation}
because $\nabla g_2(0,0)=\frac{1}{2\log 2}\big({}_{-1}^{~1}\big)$. $\square$
\section{Identification of $H(t_0)$ when $\{Y(u/2^n)\}_{u=0,\ldots,2^{n}}$ are observed}
\subsection{Estimators starting from $\{Y(u/2^n)\}_{u=0,\ldots,2^{n}}$}
Since in practice, it is more realistic to assume that a discretized trajectory of $Y(t)=\Phi(X(t))$ is observed, therefore in this section, we obtain a consistent estimator of $H(t_0)$ starting from the high frequency $\{Y(u/2^n)\}_{u=0,\ldots,2^{n}}$ is still available. Recall that the assumptions on $\Phi$ are given in Section 1.1., then the key point leading to these results is to take a second order Taylor expansion with integral remainder, to obtain for $t_0\in(0,1)$,
\begin{eqnarray}
\label{Taylor1}
&&Y(t)=\Phi(X(t_0))+\Phi'(X(t_0))(X(t)-X(t_0))\nonumber\\
&&~~+\Big(\int_0^1(1-\eta)\Phi^{(2)}(X(t_0)+\eta(X(t)-X(t_0))){\,\mathrm{d}}\eta\Big)(X(t)-X(t_0))^2.\nonumber\\
\end{eqnarray}
This together with (\ref{AH}) and the fact that  $0<c_1\leq |\Phi'(x)|\leq c_2$ for all $x\in\mathbb R$ yields
\begin{eqnarray}
\label{Taylor2}
&&\frac{d_Y(2^{-j},k)}{\Phi'(X(t_0))}=d_X(2^{-j},k)+\mathcal O_{a.s.}(2^{-j/2}\epsilon_j^{2H(t_0)}| \log\epsilon_j|).
\end{eqnarray}
Note that (\ref{AH}) can lead to
\begin{equation}
\label{bound-dX}
d_X(2^{-j},k)=\mathcal O_{a.s.}(2^{-j(H(k2^{-j})+1/2)}j^{1/2}).
\end{equation}
It follows from (\ref{Taylor2}), the definition of $V_{Y,t_0,j}$ and (\ref{bound-dX}) that
\begin{equation}
\label{VXY}
\frac{V_{Y,t_0,j}}{|\Phi'(X(t_0))|^2}=V_{X,t_0,j}\!\!+\mathcal O_{a.s.}(2^{-jH(t_0)}\epsilon_j^{2H(t_0)+1} j^{1/2}|\log\epsilon_j|+\epsilon_j^{4H(t_0)+1}|\log\epsilon_j|^2).
\end{equation}
Hence similar to Lemma \ref{IV}, we have
\begin{lemma}
\label{IVY}
If $\epsilon_j=\mathcal O(j^{-1})$, then
\begin{eqnarray}
\label{IVY1}
&&\frac{V_{Y,t_0,j}}{|\Phi'(X(t_0))|^2}=2c_2(t_0)2^{-2jH(t_0)}\epsilon_j+\mathcal O_{a.s.}\big(2^{-jH(t_0)}\epsilon_j^{2H(t_0)+1} j^{1/2}|\log\epsilon_j|\nonumber\\
&&~~~~+\epsilon_j^{4H(t_0)+1}|\log\epsilon_j|^2+2^{-2jH(t_0)}\epsilon_j(j\epsilon_j+2^{-j}j^4)\big)+\mathcal O_{\mathbb P}(2^{-j(2H(t_0)+1/2)}\epsilon_j^{1/2}).\nonumber\\
\end{eqnarray}
As a consequence,
\begin{description}
\item[(a)] Under assumption $\textbf{(A1)}$,
\begin{eqnarray}
\label{IVY2}
&&\frac{V_{Y,t_0,j}}{2|\Phi'(X(t_0))|^2c_2(t_0)2^{-2jH(t_0)}\epsilon_j}-1\nonumber\\
&&=\mathcal O_{\mathbb P}(2^{jH(t_0)}\epsilon_j^{2H(t_0)}j^{1/2}|\log\epsilon_j|+(2^{j}\epsilon_j)^{-1/2});
\end{eqnarray}
\item[(b)] Under assumptions $\textbf{(A1)-(A2)}$,
\begin{eqnarray}
\label{IVY3}
&&\frac{V_{Y,t_0,j}}{2|\Phi'(X(t_0))|^2c_2(t_0)2^{-2jH(t_0)}\epsilon_j}-1\nonumber\\
&&=\mathcal O_{a.s.}(2^{jH(t_0)}\epsilon_j^{2H(t_0)}j^{1/2}|\log\epsilon_j|+(2^{j}\epsilon_j)^{-1/2+\delta}),
\end{eqnarray}
for $\delta>0$ arbitrarily small.
\end{description}
\end{lemma}
Here we note that, to show (\ref{IVY2}) and  (\ref{IVY3}) hold, we must observe that under assumption $\textbf{(A1)}$,
$$
\epsilon_j^{4H(t_0)+1}|\log\epsilon_j|^2+2^{-2jH(t_0)}\epsilon_j(j\epsilon_j+2^{-j}j^4)=\mathcal O(2^{-jH(t_0)}\epsilon_j^{2H(t_0)+1} j^{1/2}|\log\epsilon_j|).
$$
Lemma \ref{IVY} further implies  the following:
\begin{theorem}
\label{VYd}
 For $t_0\in(0,1)$, denote by
$$
\widehat{H_{Y,2^j}}(t_0)=\frac{\log\big(\frac{\epsilon_{j+1}}{\epsilon_j}\big)+\log\big(\frac{V_{Y,t_0,j}}{V_{Y,t_0,j+1}}\big)}{2\log 2}.
$$
\begin{description}
\item[(a)] Under assumption $\textbf{(A1)}$,
\begin{equation}
\label{VYd2}
\widehat{H_{Y,2^j}}(t_0)-H(t_0)=\mathcal O_{\mathbb P}(2^{jH(t_0)}\epsilon_j^{2H(t_0)}j^{1/2}|\log\epsilon_j|+(2^{j}\epsilon_j)^{-1/2}).
\end{equation}
\item[(b)] Under assumptions $\textbf{(A1)-(A2)}$,
   \begin{equation}
\label{VYd3}
\widehat{H_{Y,2^j}}(t_0)-H(t_0)=\mathcal O_{a.s.}(2^{jH(t_0)}\epsilon_j^{2H(t_0)}j^{1/2}|\log\epsilon_j|+(2^{j}\epsilon_j)^{-1/2+\delta}),
\end{equation}
for $\delta>0$ arbitrarily small.
\item[(c)] Under assumptions $\textbf{(A2)-(A4)}$,
\begin{equation}
\label{cltY1}
\sqrt{2^{j+1}\epsilon_j}\Big(\widehat{H_{Y,2^j}}(t_0)-H(t_0)\Big)\xrightarrow[j\rightarrow+\infty]{dist}\mathcal N\big(0,\tilde{c}(t_0)\big),
\end{equation}
where the covariance matrix $\tilde{c}(t_0)$ is given in (\ref{tildec}).
\end{description}
\end{theorem}
\textbf{Proof.} (\ref{VYd2}) and (\ref{VYd3}) are obvious by using Lemma \ref{IVY} and the fact that convergences in probability and almost surely are preserved under continuous transformations. In order to show (\ref{cltY1}), we only need to verify
\begin{equation}
\label{DiffH1}
\sqrt{2^{j}\epsilon_j}|\widehat{H_{Y,2^j}}(t_0)-\widehat{H_{X,2^j}}(t_0)|\xrightarrow[j\rightarrow+\infty]{a.s.}0.
\end{equation}
Equivalently, it suffices to show
\begin{equation}
\label{DiffH2}
\sqrt{2^{j}\epsilon_j}\Big|\log\Big(\frac{V_{Y,t_0,j}}{|\Phi'(X(t_0))|^2V_{X,t_0,j}}\Big)-\log(1)\Big|\xrightarrow[j\rightarrow+\infty]{a.s.}0.
\end{equation}
We have, by using the mean value theorem on $\log(\cdot)$,
$$
\sqrt{2^{j}\epsilon_j}\Big|\log\Big(\frac{V_{Y,t_0,j}}{|\Phi'(X(t_0))|^2V_{X,t_0,j}}\Big)-\log(1)\Big|
=\frac{\sqrt{2^{j}\epsilon_j}}{\gamma_j}\Big|\frac{V_{Y,t_0,j}}{|\Phi'(X(t_0))|^2V_{X,t_0,j}}-1\Big|,
$$
where $\gamma_j$ is some random variable valued in the open interval with ending points $|\Phi'(X(t_0))|^2V_{X,t_0,j}/V_{Y,t_0,j}$ and $1$. Since $|\gamma_j|$ tends to $1$ a.s. as $j\rightarrow+\infty$, then according to (\ref{VXY}) and \textbf{(A1)}, the right-hand side of the above equation can be bounded by
$$
c2^{j(H(t_0)+1/2)}\epsilon_j^{2H(t_0)+1/2}j^{1/2}|\log\epsilon_j|~\mbox{with some}~c>0,
$$
which converges to $0$ as $j\rightarrow+\infty$, thanks to assumption \textbf{(A4)}. $\square$

The following results are of the most interests in this section. We construct consistent estimators of $H(t_0)$ starting from the observations $\{Y(u/2^n)\}_{u=0,\ldots,2^n}$.
\begin{theorem}
\label{VYdn}
 For $t_0\in(0,1)$, denote by
$$
\widehat{H_{Y,2^j,n}}(t_0)=\frac{\log\big(\frac{\epsilon_{j+1}}{\epsilon_j}\big)+\log\big(\frac{\widehat{V_{n,t_0,j}}}{\widehat{V_{n,t_0,j+1}}}\big)}{2\log 2}.
$$
\begin{description}
\item[(a)] Set $J_n=[\beta n]$, with $0<\beta<1$. We let $\epsilon_{J_n}$ satisfy assumption \textbf{(A1)}, then
\begin{eqnarray}
\label{VYdn2}
&&\widehat{H_{Y,2^{J_n},n}}(t_0)-H(t_0)=\mathcal O_{\mathbb P}(2^{(\beta-1)nH(t_0)}\nonumber\\
&&~~+2^{J_nH(t_0)}\epsilon_{J_n}^{2H(t_0)}J_n^{1/2}|\log\epsilon_{J_n}|+(2^{J_n}\epsilon_{J_n})^{-1/2}).
\end{eqnarray}
\item[(b)] Let $\epsilon_{J_n}$ satisfy assumptions \textbf{(A1)-(A2)}, then
    \begin{eqnarray}
\label{VYdn3}
&&\widehat{H_{Y,2^{J_n},n}}(t_0)-H(t_0)=\mathcal O_{a.s.}(2^{(\beta-1+\delta')nH(t_0)}\nonumber\\
&&~~+2^{J_nH(t_0)}\epsilon_{J_n}^{2H(t_0)}J_n^{1/2}|\log\epsilon_{J_n}|+(2^{J_n}\epsilon_{J_n})^{-1/2+\delta}),
\end{eqnarray}
for any $\delta,\delta'>0$ arbitrarily small.
\item[(c)] Under assumptions \textbf{(A2)-(A4)} and $0<\beta\leq\frac{4H(t_0)+1}{4H(t_0)+2}$,
\begin{equation}
\label{cltYn1}
\sqrt{2^{J_n+1}\epsilon_{J_n}}\Big(\widehat{H_{Y,2^{J_n},n}}(t_0)-H(t_0)\Big)\xrightarrow[n\rightarrow+\infty]{dist}\mathcal N\big(0,\tilde{c}(t_0)\big).
\end{equation}
\end{description}
\end{theorem}
\textbf{Proof.} In order to prove (\ref{VYdn2}) and (\ref{VYdn3}), we rely on the following relation: under assumptions \textbf{(A1)-(A2)},
\begin{equation}
\label{DiffV_p}
\frac{\widehat{V_{n,t_0,J_n}}}{V_{Y,t_0,J_n}}-1=\mathcal O_{\mathbb P}(2^{(J_n-n)H(t_0)}).
\end{equation}
This is because, by using Markov's inequality, Cauchy-Schwarz inequality, (\ref{wideV}), (\ref{IVY3}) and the dominated convergence theorem,
\begin{eqnarray}
\label{EdiffV}
&&\mathbb P\Big(\Big|\frac{\widehat{V_{n,t_0,J_n}}-V_{Y,t_0,J_n}}{V_{Y,t_0,J_n}}\Big|>\eta\Big)\leq \frac{1}{\eta}\mathbb E\Big|\frac{\widehat{V_{n,t_0,J_n}}-V_{Y,t_0,J_n}}{V_{Y,t_0,J_n}}\Big|\nonumber\\
&&\leq\frac{1}{\eta}\Big(\mathbb E\Big|\frac{\widehat{V_{n,t_0,J_n}}-V_{Y,t_0,J_n}}{2^{-2J_nH(t_0)}\epsilon_{J_n}}\Big|^2\Big)^{1/2}\Big(\mathbb E\Big|\frac{2^{-2J_nH(t_0)}\epsilon_{J_n}}{V_{Y,t_0,J_n}}\Big|^2\Big)^{1/2}\nonumber\\
&&\le c2^{(J_n-n)H(t_0)}/\eta.
\end{eqnarray}
Similarly to (\ref{EdiffV}), we also obtain for any $\delta'>0$ arbitrarily small,
$$
\mathbb P\Big(2^{(n-J_n-\delta' n)H(t_0)}\Big|\frac{\widehat{V_{n,t_0,J_n}}-V_{Y,t_0,J_n}}{V_{Y,t_0,J_n}}\Big|>\eta\Big)\leq c2^{-\delta' nH(t_0)}/\eta.
$$
The fact that $\beta<1$ implies $\sum_{n\in\mathbb N}c2^{-\delta' n H(t_0)}/\eta<+\infty$, then by Borel-Cantelli's lemma,
 \begin{equation}
 \label{DiffV_as}
 \frac{\widehat{V_{n,t_0,J_n}}}{V_{Y,t_0,J_n}}-1=\mathcal O_{a.s.}(2^{(J_n-n+\delta' n)H(t_0)})=\mathcal O_{a.s.}(2^{(\beta-1+\delta')nH(t_0)}).
 \end{equation}
  Therefore, (\ref{VYdn2}) (resp. (\ref{VYdn3})) follows from the following 2 decompositions:
  $$
  \widehat{H_{Y,2^{J_n},n}}(t_0)-H(t_0)=\widehat{H_{Y,2^{J_n},n}}(t_0)-\widehat{H_{Y,2^{J_n}}}(t_0)+\widehat{H_{Y,2^{J_n}}}(t_0)-H(t_0);
  $$
\begin{eqnarray*}
&&\frac{\widehat{V_{n,t_0,J_n}}}{2|\Phi'(X(t_0))|^2c_2(t_0)2^{-2J_nH(t_0)}\epsilon_{J_n}}-1\nonumber\\
&&=
\Big(\frac{V_{Y,t_0,J_n}}{2|\Phi'(X(t_0))|^2c_2(t_0)2^{-2J_nH(t_0)}\epsilon_j}-1\Big)\Big(\frac{\widehat{V_{n,t_0,J_n}}}{V_{Y,t_0,J_n}}\Big)
+\Big(\frac{\widehat{V_{n,t_0,J_n}}-V_{Y,t_0,J_n}}{V_{Y,t_0,J_n}}\Big);
\end{eqnarray*}
and equations (\ref{VYd2}), (\ref{DiffV_p}) (resp. (\ref{VYd3}), (\ref{DiffV_as})).

For showing (\ref{cltYn1}), we only need to show
\begin{equation}
\label{DiffH}
\sqrt{2^{J_n+1}\epsilon_{J_n}}\big|\widehat{H_{Y,2^{J_n},n}}(t_0)-\widehat{H_{Y,2^{J_n}}}(t_0)\big|\xrightarrow[n\rightarrow+\infty]{\mathbb P}0.
\end{equation}
By using the same idea we took to prove (\ref{DiffH1}), we just need to verify
$$
\sqrt{2^{J_n+1}\epsilon_{J_n}}\Big|\frac{\widehat{V_{n,t_0,J_n}}-V_{Y,t_0,J_n}}{V_{Y,t_0,J_n}}\Big|\xrightarrow[n\rightarrow+\infty]{\mathbb P}0.
$$
This is true since, according to (\ref{EdiffV}) and the fact that (by assumption \textbf{(A4)}) $\epsilon_{J_n}=o\big(2^{-(\frac{2H(t_0)+1}{4H(t_0)+1})J_n}\big)$,
the left-hand side of the above term can be bounded in probability by
$$
c2^{J_n/2}\epsilon_{J_n}^{1/2}2^{(J_n-n)H(t_0)}=o\big(2^{H(t_0)(J_n\frac{4H(t_0)+2}{4H(t_0)+1}-n)}\big).
$$
The fact that $0<\beta\leq\frac{4H(t_0)+1}{4H(t_0)+2}$ entails $J_n\frac{4H(t_0)+2}{4H(t_0)+1}-n\leq 0$.  Consequently,
$$
\sqrt{2^{J_n+1}\epsilon_{J_n}}\Big|\frac{\widehat{V_{n,t_0,J_n}}-V_{Y,t_0,J_n}}{V_{Y,t_0,J_n}}\Big|\xrightarrow[n\rightarrow+\infty]{\mathbb P}0,
$$
hence (\ref{cltYn1}) holds. $\square$
\subsection{Selection of the parameters $\epsilon_j$ and $\beta$ from practical point of view}
In Theorem \ref{VYdn}, the choices of $\epsilon_j$ and $\beta$ depend on the target parameter $H(t_0)$. This is unacceptable from practical point of view. To overcome this inconvenience, we make assumptions \textbf{(A1), (A4)} stronger so that the values of $\epsilon_j$ and $\beta$ don't rely on $H(t_0)$: suppose the lower bound $\tau_1=\inf\limits_{t\in[0,1]}H(t)$ and the upper bound $\tau_2=\sup\limits_{t\in[0,1]}H(t)$ are known.
 \begin{description}
 \item[\textbf{(A1)\textprime}] $\epsilon_j=2^{-j\gamma}$ with $\gamma\in(\frac{1}{2},1)$;
 \item[\textbf{(A4)\textprime}] $\epsilon_j=2^{-j\gamma}$ with $\gamma\in(\frac{2\tau_2+1}{4\tau_2+1},1)$ and $0<\beta\leq\frac{4\tau_1+1}{4\tau_1+2}$.
 \end{description}
 Without great effort we could see \textbf{(A1)\textprime} implies \textbf{(A1)} and \textbf{(A4)\textprime} implies \textbf{(A4)} for all $t_0\in(0,1)$. Then from Theorem \ref{VYd} and Theorem \ref{VYdn}, we easily derive the following results:
 \begin{corollary}
 \label{cor2}
 \begin{description}
\item[(a)] Let $\epsilon_{j}$ satisfy assumption \textbf{(A1)\textprime}, then
$$
\sup_{t\in(0,1)}|\widehat{H_{Y,2^j}}(t)-H(t)|=\mathcal O_{\mathbb P}\big(2^{j(1-2\gamma)\tau_1}j^{3/2}+2^{-j(1-\gamma)/2}\big).
$$
\item[(b)] If $\epsilon_{j}$ satisfies assumption \textbf{(A1)\textprime}, then it also satisfies \textbf{(A2)}. As a result,
  $$
  \sup_{t\in(0,1)}|\widehat{H_{Y,2^j}}(t)-H(t)|=\mathcal O_{a.s.}\big(2^{j(1-2\gamma)\tau_1}j^{3/2}+2^{j(1-\gamma)(-1/2+\delta)}\big)
  $$
for $\delta>0$ arbitrarily small.
\item[(c)] Under assumptions \textbf{(A2), (A3), (A4)\textprime},
$$
\sqrt{2^{j+1}\epsilon_j}\Big(\widehat{H_{Y,2^j}}(t_0)-H(t_0)\Big)\xrightarrow[j\rightarrow+\infty]{dist}\mathcal N\big(0,\tilde{c}(t_0)\big).
$$
\end{description}
 \end{corollary}
 \begin{corollary}
 \label{cor1}
 \begin{description}
\item[(a)] Set $J_n=[\beta n]$, with $0<\beta<1$. If $\epsilon_{J_n}$ satisfies assumption \textbf{(A1)\textprime}, then
$$
\sup_{t\in(0,1)}|\widehat{H_{Y,2^{J_n},n}}(t)-H(t)|=\mathcal O_{\mathbb P}(2^{(\beta-1)n\tau_1}+2^{\beta(1-2\gamma)n\tau_1}n^{3/2}+2^{\beta(\gamma-1)n/2}).
$$
\item[(b)] If $\epsilon_{J_n}$ satisfies assumption \textbf{(A1)\textprime}, then
   \begin{eqnarray*}
&&\sup_{t\in(0,1)}|\widehat{H_{Y,2^{J_n},n}}(t)-H(t)|\\
&&=\mathcal O_{a.s.}(2^{(\beta-1+\delta')n\tau_1}+2^{\beta(1-2\gamma)n\tau_1}n^{3/2}+2^{\beta(\gamma-1)n(1/2-\delta)}),
\end{eqnarray*}
for any $\delta,\delta'>0$ arbitrarily small.
\item[(c)] Under assumptions \textbf{(A2), (A3), (A4)\textprime},
$$
\sqrt{2^{J_n+1}\epsilon_{J_n}}\Big(\widehat{H_{Y,2^{J_n},n}}(t_0)-H(t_0)\Big)\xrightarrow[n\rightarrow+\infty]{dist}\mathcal N\big(0,\tilde{c}(t_0)\big).
$$
\end{description}
 \end{corollary}
We conclude that the selection of $\beta,\gamma$ can be made for each setting as follows:
 \begin{itemize}
 \item In Corollary \ref{cor2} (a), for any $\tau_1\in(0,1)$, taking $\gamma=\frac{2\tau_1+1}{4\tau_1+1}$, we obtain the best rate of convergence in this setting:
   $$
\sup_{t\in(0,1)}|\widehat{H_{Y,2^j}}(t)-H(t)|=\mathcal O_{\mathbb P}\big(2^{-j(\frac{\tau_1}{4\tau_1+1})}j^{3/2}\big).
$$
\item In Corollary \ref{cor2} (b), similarly, taking $\gamma=\frac{2\tau_1+1-2\delta}{4\tau_1+1-2\delta}$ obtains
   $$
\sup_{t\in(0,1)}|\widehat{H_{Y,2^j}}(t)-H(t)|=\mathcal O_{a.s.}\big(2^{-j(\frac{\tau_1(1+2\delta)}{4\tau_1+1-2\delta})}j^{3/2}\big).
$$
   \item In Corollary \ref{cor1} (a), the best choices of $\beta$, $\gamma$ are the ones such that $\beta=\frac{4\tau_1+1}{4\tau_1+2}$ and $\gamma=\frac{2\tau_1+1}{4\tau_1+1}$. Consequently,
   $$
\sup_{t\in(0,1)}|\widehat{H_{Y,2^{J_n},n}}(t)-H(t)|=\mathcal O_{\mathbb P}(2^{-(\frac{\tau_1}{4\tau_1+2})n}n^{3/2}).
$$
\item In Corollary \ref{cor1} (b) we select $\beta=\frac{(1-\delta')(4\tau_1+1-2\delta)}{4\tau_1+2-4\delta}$ and $\gamma=\frac{2\tau_1+1-2\delta}{4\tau_1+1-2\delta}$ to obtain
   $$
\sup_{t\in(0,1)}|\widehat{H_{Y,2^{J_n},n}}(t)-H(t)|=\mathcal O_{a.s.}(2^{-(\frac{(1-\delta')(1-2\delta)\tau_1}{4\tau_1+2-4\delta})n}n^{3/2}).
$$
 \end{itemize}
 From the above discussion we conclude that our estimator has its best convergence rate approximately $2^{-j(\frac{\tau_1}{4\tau_1+1})}$ when the wavelet coefficients of $Y$ are observed, and $2^{-n(\frac{\tau_1}{4\tau_1+2})}$ when a discrete sample path of $Y$ is observed. This convergence rate is poorer than the case when $X$ is straightforwardly available.
\section{An application to statistical inferences: estimation of $\theta(\cdot)$}
The estimation of the hidden pointwise H\"older exponent also allows to solve stochastic volatility model's nonparametric estimation problem. For example, suppose some underlying mBm $B_{H(t)}(t)$ satisfies the following nonlinear model: for $t\in[0,1]$,
\begin{equation}
\label{model}
Y(t)=\theta(t) B_{H(t)}(t),
\end{equation}
where
\begin{itemize}
  \item $\{B_{H(t)}(t)\}_{t\in[0,1]}$ is an mBm with unknown $H\in C^2([0,1])$;
  \item $\theta\in C^2([0,1])$, $\theta\neq0$ almost everywhere is an unknown real-valued deterministic function;
  \item only a discrete trajectory $\{Y(u/2^n)\}_{u=0,\ldots,2^n}$ is available.
\end{itemize}
From the statistical setting, the information of the hidden mBm remains unknown. However, the following result shows it is still possible to evaluate the parameter $\theta$. The strategy is: first estimate the pointwise H\"older exponent $H(\cdot)$ by $\widehat{H_{Y,2^{J_n},n}}(\cdot)$ given in Theorem \ref{VYdn}, then the following result holds:
\begin{proposition}
\label{prop:application}
Fix $t_0\in(0,1)$ and $J_n=[\beta n]$ with $\beta\in(0,1)$. Let
$$
\widehat{\theta_n}^2(t_0)=-\frac{4^{J_n\widehat{H_{Y,2^{J_n},n}}(t_0)}\epsilon_{J_n}^{-1}\widehat{V_{n,t_0,J_n}}}{C_0(2\widehat{H_{Y,2^{J_n},n}}(t_0))\int_0^1\int_0^1
 \psi(t)\psi(s)|t-s|^{2\widehat{H_{Y,2^{J_n},n}}(t_0)}{\,\mathrm{d}} t{\,\mathrm{d}} s}.
$$
 Then,
 \begin{description}
 \item[(a)] Under assumption $\textbf{(A1)\textprime}$,
\begin{equation}
\label{theta1'}
\widehat{\theta_n}(t_0)^2\xrightarrow[n\rightarrow+\infty]{\mathbb P}\theta(t_0)^2.
\end{equation}
\item[(b)]  Under assumptions $\textbf{(A1)\textprime}$ and \textbf{(A2)},
\begin{equation}
\label{theta1}
\widehat{\theta_n}(t_0)^2\xrightarrow[n\rightarrow+\infty]{a.s.}\theta(t_0)^2.
\end{equation}
\end{description}
\end{proposition}
\textbf{Proof of Proposition \ref{prop:application}.} We only show (\ref{theta1}) holds, since the way to obtain (\ref{theta1'}) is similar. Under assumptions \textbf{(A1)\textprime} and \textbf{(A2)}, on one hand, it follows from (\ref{DiffV_as}) and (\ref{IVY3}) that
\begin{equation}
\label{ptheta1}
\frac{\widehat{V_{n,t_0,J_n}}}{2^{-2J_nH(t_0)}\epsilon_{J_n}}\xrightarrow[J_n\rightarrow+\infty]{a.s.}2c_2(t_0).
\end{equation}
On the other hand, by using (\ref{c2}) and the fact that $c_2(t_0)=C_2(t_0,1)$, we see
\begin{equation}
\label{diffhatH}
\theta(t_0)^2=-\frac{2c_2(t_0)}{C_{0}(2H(t_0))\int_{0}^{1}\int_{0}^{1}\psi
( t) \psi ( s) | t-s| ^{2H(t_0
) }{\,\mathrm{d}} t{\,\mathrm{d}} s}.
\end{equation}
Since the functions $x\mapsto C_0(2x)$ and $x\mapsto \int_{0}^{1}\int_{0}^{1}\psi( t) \psi ( s) | t-s| ^{2x}{\,\mathrm{d}} t{\,\mathrm{d}} s$
are both continuous over $(0,1)$, therefore by combining (\ref{ptheta1}), (\ref{diffhatH}) and (\ref{VYdn3}), we obtain Proposition \ref{prop:application}. $\square$
\section{Comparison to Bardet and Surgailis (2013) and a simulation study}
\label{simulation}
Bardet and Surgailis (2013) established a pseudo-generalized least squares version of the localized increment ratio estimator and quadratic estimator of $H(\cdot)$, denoted by $\widehat{H}_{n,\alpha}^{IR2}$ and $\widehat{H}_{n,\alpha}^{QV2}$ respectively. This work is so far the most achieved one on estimation of the mBm's pointwise H\"older exponent. In this section we compare our model and approach to Bardet and Surgailis (2013)'s and illustrate some simulation results. The main differences between our statistical setting and Bardet and Surgailis (2013)'s are:
\begin{enumerate}
\item Bardet and Surgailis (2013) considers observation of a multifractional Gaussian process which has asymptotic self-similarity and \enquote{tangent to} an fBm with Hurst parameter $H(t)$ and scaled by $c(t)$ at each time $t$. This setting covers ours when $\Phi=Id$. However, our model is different when $\Phi$ is some other function.
\item The estimators $\widehat{H}_{n,\alpha}^{IR2}$ and $\widehat{H}_{n,\alpha}^{QV2}$ are obtained in terms of observations of discrete sample path of the process. We have established estimators based on observations of discrete sample path of $Y$ and wavelet coefficients of the process (see $\widehat{H_{Y,2^j}}$ in Theorem \ref{VYd}), respectively. The latter result is the main contribution of this paper. For more details on wavelet-based statistics we refer to Delbeke and Van (1995), Abry et al. (2002), Abry and Con\c{c}alv\`es (1997).
\item The estimators $\widehat{H}_{n,\alpha}^{IR2}$ and $\widehat{H}_{n,\alpha}^{QV2}$ apply to all $H\in C^{\eta}([0,1])$ with $\eta>0$, however our estimators only work for $\eta\ge2$.
    \item In both $\widehat{H}_{n,\alpha}^{IR2}$ and $\widehat{H}_{n,\alpha}^{QV2}$, the length of the local window $n^{1-\alpha}$ plays the role as $2^j\epsilon_j$ in our setting, which partially explains the estimator's asymptotic behavior. For building up the convergence, both Bardet and Surgailis (2013)'s and our model are subject to some technical assumptions. Here we consider a particular case to compare the estimators' rate of convergence. Let $Y(t)=\theta(t)X(t)$, where $X$ has pointwise H\"older exponent $H(\cdot)$. From Proposition 3 (iii) and (iv) in Bardet and Surgailis (2013) we see for $\eta=2$ and arbitrarily small $\epsilon,\delta>0$, the estimators $\widehat H_{n,2/11}^{\mbox {\textbf{$\varepsilon$}}}(t)$ and $\widehat H_{n,1/11}^{\mbox {\textbf{$\varepsilon$}}}(t)$ have the best rate of convergence in probability and a.s. respectively:
        $$
        \sup_{\epsilon<t<1-\epsilon}|\widehat H_{n,2/11}^{\mbox {\textbf{$\varepsilon$}}}(t)-H(t)|=\mathcal O_{\mathbb P}(n^{-4/11})$$
        and
        $$
        \sup_{\epsilon<t<1-\epsilon}|\widehat H_{n,1/11}^{\mbox {\textbf{$\varepsilon$}}}(t)-H(t)|=\mathcal O_{a.s.}(n^{-2/11+\delta}),
        $$
        where \textbf{$\varepsilon$}$=$ IR2 or QV2. We set $j=\log n/\log 2$ in Corollary \ref{cor}. Then according to Corollary \ref{cor}, $\widehat H_{n,2/11}^{\mbox {\textbf{$\varepsilon$}}}(t)$ has a better rate of convergence in probability than $\widehat{H_{X,2^j}}(t)$, while $\widehat{H_{X,2^j}}(t)$ outperforms $\widehat H_{n,1/11}^{\mbox {\textbf{$\varepsilon$}}}(t)$ by having a better rate of a.s. convergence.
\end{enumerate}
Now we deliver algorithms to generate our estimators $\widehat{H_{Y,2^j}}(t_0)$ and $\widehat{H_{Y,t_0,n}}(t_0)$. From Theorem \ref{VYd}, an algorithm to generate $\widehat{H_{Y,2^j}}(t_0)$ can be given as follows:\\
\noindent\makebox[\linewidth]{\rule{\textwidth}{0.4pt}}
\textbf{Algorithm I: estimation of $ H(t_0)$ starting from wavelet coefficients }\\
\noindent\makebox[\linewidth]{\rule{\textwidth}{0.4pt}}
$1$: \textbf{INPUT:} a sample $\{d_Y(2^{-j},k)\}_{k=0,\ldots,2^{j}-1}$, $\tau_1=\inf\limits_{t\in[0,1]}H(t)$ and $t_0\in(0,1)$\\
$2$: \emph{\textbf{Initialize all the parameters:}}\\
$3$: $\gamma\leftarrow\frac{2\tau_1+1}{4\tau_1+1}$; $\epsilon_{j}\leftarrow2^{-j\gamma}$;\\
$4$: \emph{\textbf{Establish a set of indices corresponding to the neighbors of $Y(t_0)$:}}\\
$5$: $\nu_{t_0,2^{j}}\leftarrow\big\{0,\ldots,2^{j}-1\big\}\cap\big[2^{j}(t_0-\epsilon_{j}),2^{j}(t_0+\epsilon_{j})\big]$;\\
$6$: \emph{\textbf{Compute \enquote{partial} sum of squares of the wavelet coefficients:}}\\
$7$: $V_{Y,t_0,j}\leftarrow\sum\limits_{k\in\nu_{t_0,2^j}}d_Y(2^{-j},k)^2$;\\
$8$: \emph{\textbf{Compute estimate of $H(t_0)$:}}\\
$9$: $\widehat{H_{Y,2^j}}(t_0)\leftarrow\frac{1}{2\log(2)}\log\big(\frac{\epsilon_{j+1}{V_{Y,t_0,j}}}{\epsilon_{j}{V_{Y,t_0,j+1}}}\big)$;\\
$10$: \textbf{OUTPUT:} $\widehat{H_{Y,2^j}}(t_0)$\\
\noindent\makebox[\linewidth]{\rule{\textwidth}{0.4pt}}
In Algorithm I, Line 3 proposes a data-driven procedure to set the parameters $\gamma$ and  $\varepsilon_j$. From Line 7 we see that this algorithm requires at least $\mathcal O(\nu_{t_0,2^{j}})=\mathcal O(2^{j(\frac{2\tau_1}{4\tau_1+1})})$ operations. By this result and Corollary \ref{cor2} it appears a trade-off between computational cost and precision of estimation: less value of $\tau_1$ reduces computational cost but leads to less precision of estimation.

Next we numerically study the estimator $\widehat{H_{Y,2^{J_n},n}}(t_0)$ and compare the results with Bardet and Surgailis (2013)'s. In view of Corollary  \ref{cor1}, we choose the parameters $J_n=[\beta n]$, $\epsilon_{J_n}=2^{-J_n\gamma}$ with $\beta=\frac{4\tau_1+1}{4\tau_1+2}$ and $\gamma=\frac{1}{2\beta}$. We also choose $\psi$ to be Haar wavelet's mother wavelet function: $\psi(x)=\mathds{1}_{[0,1/2)}(x)-\mathds{1}_{[1/2,1]}(x)$. Below is an algorithm to simulate $\widehat{H_{Y,2^{J_n},n}}(t_0)$, according to  Theorem \ref{VYdn}.\\

\noindent\makebox[\linewidth]{\rule{\textwidth}{0.4pt}}
\textbf{Algorithm II: estimation of $H(t_0)$ starting from discretized trajectory}\\
\noindent\makebox[\linewidth]{\rule{\textwidth}{0.4pt}}
$1$: \textbf{INPUT:} a sample $\{Y(k/2^n)\}_{k=0,\ldots,2^n}$, $\tau_1=\inf\limits_{t\in[0,1]}H(t)$ and $t_0\in(0,1)$\\
$2$: \emph{\textbf{Initialize all the parameters:}}\\
$3$: $\beta\leftarrow\frac{4\tau_1+1}{4\tau_1+2}$; $\gamma\leftarrow\frac{1}{2\beta}$; $J_n\leftarrow[\beta n]$; $\epsilon_{J_n}\leftarrow2^{-J_n\gamma}$;\\
$4$: \emph{\textbf{Establish a set of indices corresponding to the neighbors of $Y(t_0)$:}}\\
$5$: $\nu_{t_0,2^{J_n}}\leftarrow\big\{0,\ldots,2^{J_n}-1\big\}\cap\big[2^{J_n}(t_0-\epsilon_{J_n}),2^{J_n}(t_0+\epsilon_{J_n})\big]$;\\
$6$: \emph{\textbf{Estimate the wavelet coefficients:}}\\
$7$: \textbf{for $k$ in $\nu_{t_0,2^{J_n}}$, do:}\\
$8$: ~~~~~~~~$\widehat{d_{Y,n}}(2^{-J_n},k)\leftarrow2^{J_n/2}\sum\limits_{l=0}^{2^{n-J_n}-1}Y(l2^{-n}+k2^{-J_n})\int_{\frac{l}{2^n}}^{\frac{l+1}{2^n}}\psi(2^{J_n}t){\,\mathrm{d}} t$;\\
$9$: \textbf{end for}\\
$10$: \emph{\textbf{Estimate \enquote{partial} sum of squares of the wavelet coefficients:}}\\
$11$: $ \widehat{V_{n,t_0,J_n}}\leftarrow\sum\limits_{k\in\nu_{t_0,2^{J_n}}}\widehat{d_{Y,n}}(2^{-J_n},k)^2$;\\
$12$: \emph{\textbf{Compute estimate of $H(t_0)$:}}\\
$13$: $\widehat{H_{Y,2^{J_n},n}}(t_0)\leftarrow\frac{1}{2\log(2)}\log\big(\frac{\epsilon_{J_n+1}\widehat{V_{n,t_0,J_n}}}{\epsilon_{J_n}\widehat{V_{n,t_0,J_n+1}}}\big)$;\\
$14$: \textbf{OUTPUT:} $\widehat{H_{Y,2^{J_n},n}}(t_0)$\\
\noindent\makebox[\linewidth]{\rule{\textwidth}{0.4pt}}
Line 3 shows the procedure to set the parameters $\beta$, $\gamma$, $J_n$ and $\epsilon_{J_n}$. In Line 8, each item $\int_{\frac{l}{2^n}}^{\frac{l+1}{2^n}}\psi(2^{J_n}t){\,\mathrm{d}}t$ can be approximated by $\psi(2^{J_n-n}(l+1/2))$, when $n$ is large. Hence from Lines 7-9 of the algorithm we see that the algorithm requires at least $\mathcal O(\nu_{t_0,2^{J_n}}2^{n-J_n})=\mathcal O(2^{n/2})$ operations. This fact together with the discussion on Corollary \ref{cor1} reveals that the choice of $\tau_1$ has no impact on the computational cost, but only on the precision of the estimation.

Next we provide the empirical mean, the standard deviation and the quantile-quantile plots (QQ plots) through a simulation study, where all the codes in MATLAB are available from authors upon request.

 We let $\theta(\cdot)\equiv1$, then choose 4 different types of $\Phi$ and $3$ different Hurst functional parameters $H(t)$, $t\in(0,1)$ as follows:
 \begin{eqnarray*}
 &&\Phi_1(x)=x;~\Phi_2(x)=e^x;~\Phi_3(x)=\sin(4x);~\Phi_4(x)=x\sin^2(4x);\\
 &&H_1(t)=0.1+0.8t;~H_2(t)=0.5+0.4\sin(5t);~H_3(t)=0.1+0.8(1-t)\sin^2(10t).
 \end{eqnarray*}
 Since $\tau_1=0.1$ for the above three $H(\cdot)$, we then choose $Q=2$, $\beta=\frac{4\tau_1+1}{4\tau_1+2}=\frac{7}{12}$,  $\gamma=\frac{1}{2\beta}=\frac{6}{7}$. We assume a single discrete trajectory:  $(Y(k2^{-13}))_{k=0,\ldots,2^{13}}$ is available, and denote $\widehat{H_i}(t_0)$ to be the estimate of $H_i(t_0)$: $\widehat{H_{Y,2^{J_n},n}^{(i)}}(t_0)$ ($i=1,2,3$) for short. By generating each estimator $100$ times, we present the corresponding empirical mean ($m$) and standard deviation ($std$) in the following table.
\begin{table}[H]\setlength\extrarowheight{1pt}
  \begin{tabular}{|l|l|l|l|l|l|l|}
    \hline
    \multirow{2}{1.5cm}{$\Phi_1$} &
      \multicolumn{2}{c|}{$\widehat{H_{1}}(t_0)$} &
      \multicolumn{2}{c|}{$\widehat{H_{2}}(t_0)$} &
      \multicolumn{2}{c|}{$\widehat{H_{3}}(t_0)$} \\
      \cline{2-7}
    & $m\backslash H_1(t_0)$ & $std$ & $m\backslash H_2(t_0)$ & $std$ & $m\backslash H_3(t_0)$ & $std$ \\
    \hline
    $t_0=0.1$ & $0.1501\backslash0.18$ & 0.0636 & $0.6638\backslash0.6918$ & 0.0851 & $0.6362\backslash0.6098$ & 0.0874 \\
    \hline
    $t_0=0.3$ & $0.3253\backslash0.34$ & 0.0612 & $0.8860\backslash0.8990$ & 0.0487 & $0.1134\backslash0.1112$ & 0.0618 \\
    \hline
    $t_0=0.5$ & $0.4919\backslash0.50$ & 0.0775 & $0.7178\backslash0.7394$ & 0.0699 & $0.4455\backslash0.4678$ & 0.0676 \\
    \hline
    $t_0=0.7$ & $0.6413\backslash0.66$ & 0.0787 & $0.3785\backslash0.3597$ & 0.0696 & $0.1897\backslash0.2036$ & 0.0491 \\
    \hline
    $t_0=0.9$ & $0.7972\backslash0.82$ & 0.0939 & $0.0714\backslash0.1090$ & 0.0609 & $0.0848\backslash0.1136$ & 0.0605 \\
    \hline
   \hline
    \multirow{2}{1cm}{$\Phi_2$} &
      \multicolumn{2}{c|}{$\widehat{H_{1}}(t_0)$} &
      \multicolumn{2}{c|}{$\widehat{H_{2}}(t_0)$} &
      \multicolumn{2}{c|}{$\widehat{H_{3}}(t_0)$} \\
      \cline{2-7}
    & $m\backslash H_1(t_0)$ & $std$ & $m\backslash H_2(t_0)$ & $std$ & $m\backslash H_3(t_0)$ & $std$ \\
    \hline
    $t_0=0.1$ & $0.1560\backslash0.18$ & 0.0769 & $0.6518\backslash0.6918$ & 0.0780 & $0.6371\backslash0.6098$ & 0.0779 \\
    \hline
    $t_0=0.3$ & $0.3395\backslash0.34$ & 0.0683 & $0.8887\backslash0.8990$ & 0.0568 & $0.1104\backslash0.1112$ & 0.0845 \\
    \hline
    $t_0=0.5$ & $0.5139\backslash0.50$ & 0.0711 & $0.7234\backslash0.7394$ & 0.0897 & $0.4528\backslash0.4678$ & 0.0805 \\
    \hline
    $t_0=0.7$ & $0.6228\backslash0.66$ & 0.0823 & $0.3705\backslash0.3597$ & 0.0874 & $0.2060\backslash0.2036$ & 0.0767 \\
    \hline
    $t_0=0.9$ & $0.7712\backslash0.82$ & 0.0609 & $0.0667\backslash0.1090$ & 0.0666 & $0.0989\backslash0.1136$ & 0.0774 \\
    \hline
   \hline
    \multirow{2}{1cm}{$\Phi_3$} &
      \multicolumn{2}{c|}{$\widehat{H_{1}}(t_0)$} &
      \multicolumn{2}{c|}{$\widehat{H_{2}}(t_0)$} &
      \multicolumn{2}{c|}{$\widehat{H_{3}}(t_0)$} \\
      \cline{2-7}
    & $m\backslash H_1(t_0)$ & $std$ & $m\backslash H_2(t_0)$ & $std$ & $m\backslash H_3(t_0)$ & $std$ \\
    \hline
    $t_0=0.1$ & $0.1548\backslash0.18$ & 0.0695& $0.7420\backslash0.6918$ & 0.1647 & $0.6588\backslash0.6098$ &0.1185  \\
    \hline
    $t_0=0.3$ & $0.2935\backslash0.34$ & 0.0835 & $0.8395\backslash0.8990$ & 0.1043 & $0.1138\backslash0.1112$ &  0.0580 \\
    \hline
    $t_0=0.5$ & $0.5389\backslash0.50$ & 0.1150 & $0.7857\backslash0.7394$ & 0.1064 & $0.4449\backslash0.4678$ & 0.1075\\
    \hline
    $t_0=0.7$ & $0.6902\backslash0.66$ & 0.1200 & $0.3176\backslash0.3597$ & 0.0971 & $0.2234\backslash0.2036$ & 0.0681 \\
    \hline
    $t_0=0.9$ & $0.7692\backslash0.82$ & 0.1092 & $0.1083\backslash0.1090$ & 0.0563 & $0.0907\backslash0.1136$ & 0.0628\\
    \hline
     \hline
    \multirow{2}{1cm}{$\Phi_4$} &
      \multicolumn{2}{c|}{$\widehat{H_{1}}(t_0)$} &
      \multicolumn{2}{c|}{$\widehat{H_{2}}(t_0)$} &
      \multicolumn{2}{c|}{$\widehat{H_{3}}(t_0)$} \\
      \cline{2-7}
    & $m\backslash H_1(t_0)$ & $std$ & $m\backslash H_2(t_0)$ & $std$ & $m\backslash H_3(t_0)$ & $std$ \\
    \hline
    $t_0=0.1$ & $0.1589\backslash0.18$ & 0.0713& $0.6812\backslash0.6918$ & 0.0852 & $0.6924\backslash0.6098$ & 0.1196 \\
    \hline
    $t_0=0.3$ & $0.3214\backslash0.34$ & 0.0790 & $0.8505\backslash0.8990$ & 0.0564 & $0.0843\backslash0.1112$ &0.0686 \\
    \hline
    $t_0=0.5$ & $0.5280\backslash0.50$ & 0.0988 & $0.7459\backslash0.7394$ & 0.1072 & $0.5366\backslash0.4678$ & 0.1267  \\
    \hline
    $t_0=0.7$ & $0.6682\backslash0.66$ & 0.0938 & $0.3704\backslash0.3597$ & 0.0908 & $0.1929\backslash0.2036$ &0.0579  \\
    \hline
    $t_0=0.9$ & $0.7478\backslash0.82$ & 0.1145 & $0.1045\backslash0.1090$ & 0.0610 & $0.0953\backslash0.1136$ & 0.0619  \\
    \hline
  \end{tabular}
   \caption{Mean and std of the estimates $\widehat{H_{Y,2^{J_n},n}}(t_0)$ with $n=13$ and $t_0=0.1,0.3,0.5,0.7,0.9$.}
\end{table}
\noindent
From Table 1 we see that our estimator has a consistent performance on evaluating the true pointwise H\"older exponent with respect to different $H$ and $\Phi$.
Now we consider comparing the IR2 estimators in Bardet and Surgailis (2013) with ours. To this end we take $\Phi_1(x)=x$ to let $Y$ be an mBm. By using the MATLAB functions \emph{VariaIR$\_$MBM(eta,$n$,tot)} and \emph{HestIR($R$,$p$)} obtained from  \emph{http://samm.univ-paris1.fr/-Jean-Marc-Bardet}, we also provide the empirical mean and standard deviation (based on $100$ independent scenarios of mBm and $\alpha=0.5$) of IR2 estimators as follows:
\begin{table}[H]\setlength\extrarowheight{1pt}
  \begin{tabular}{|l|l|l|l|l|l|l|}
    \hline
    \multirow{2}{1.5cm}{mBm} &
      \multicolumn{2}{c|}{$\widehat{H_{1}^{IR2}}(t_0)$} &
      \multicolumn{2}{c|}{$\widehat{H_{2}^{IR2}}(t_0)$} &
      \multicolumn{2}{c|}{$\widehat{H_{3}^{IR2}}(t_0)$} \\
      \cline{2-7}
    & $m\backslash H_1(t_0)$ & $std$ & $m\backslash H_2(t_0)$ & $std$ & $m\backslash H_3(t_0)$ & $std$ \\
    \hline
    $t_0=0.1$ & $0.2070\backslash0.18$ & 0.0513 & $0.8094\backslash0.6918$ & 0.0801 & $0.6661\backslash0.6098$ & 0.0636 \\
    \hline
    $t_0=0.3$ & $0.3356\backslash0.34$ & 0.0557 & $0.8581\backslash0.8990$ & 0.0838 & $0.1291\backslash0.1112$ & 0.0450 \\
    \hline
    $t_0=0.5$ & $0.5308\backslash0.50$ & 0.0652 & $0.7163\backslash0.7394$ & 0.0810 & $0.4498\backslash0.4678$ & 0.0759 \\
    \hline
    $t_0=0.7$ & $0.6901\backslash0.66$ & 0.0664 & $0.3493\backslash0.3597$ & 0.0606 & $0.1909\backslash0.2036$ & 0.0603 \\
    \hline
    $t_0=0.9$ & $0.7778\backslash0.82$ & 0.0794 & $0.1047\backslash0.1090$ & 0.0435 & $0.1330\backslash0.1136$ & 0.0509 \\
    \hline
  \end{tabular}
   \caption{Mean and std of the IR2 estimators with $n=2^{13}$ and $t_0=0.1,0.3,0.5,0.7,0.9$.}
\end{table}
\noindent
Note that in the function \emph{VariaIR$\_$MBM(eta,n,tot)}, the mBm was previously generated using  the Choleski decomposition of the covariance matrix, so that the sample size $n$ is limited to $6000$. However here we have generated the mBm using  Wood $\&$ Chan circulant matrix, some krigging and a prequantification (see Chan and Wood 1998; Barri\`ere 2007), the sample size $n$ can be thus taken as $2^{13}=8192$. The empirical comparison shows no significant difference between the performances of $\widehat{H_{Y,2^{J_n},n}}(t_0)$ and $\widehat{H_{n,\alpha}^{IR2}}(t_0)$, except that $IR2$ estimator has less variance. Below we compare the probability distributions of $\widehat{H_{Y,2^{J_n},n}}(t_0)$ versus $\widehat{H_{n,\alpha}^{IR2}}(t_0)$, by displaying QQ plots (see Fig. 1). The first $3$ QQ plots show whether $\widehat{H_{Y,2^{J_n},n}}(t_0)$ and $\widehat{H_{n,\alpha}^{IR2}}(t_0)$ come from the same distribution for $n=2^{13}$, $\alpha=0.5$, $H(t)=0.1+0.8t$ and $t_0=0.3,0.5,0.7$; the last one illustrates the asymptotically normal behavior of our estimator $\widehat{H_{Y,2^{J_n},n}}(t_0)$ for $H(t)=0.1+0.8t$ and $t_0=0.5$.

\begin{figure}[H]
\includegraphics[scale=.4]{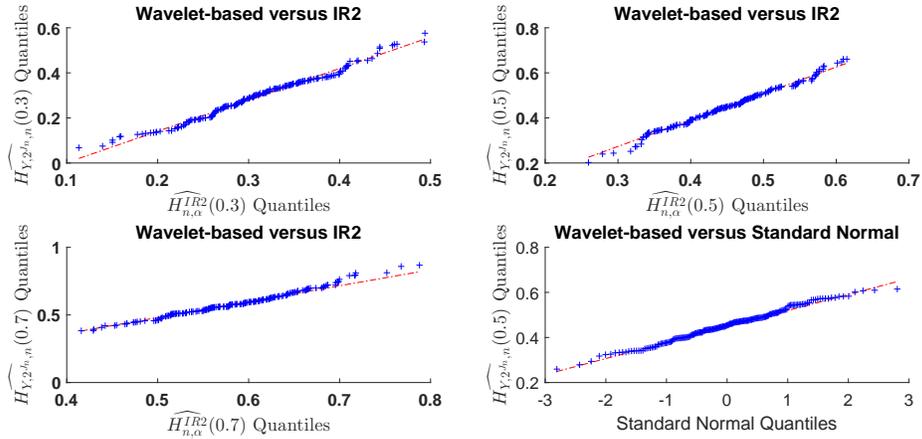}
\caption{QQ plot of $\widehat{H_{Y,2^{J_n},n}}(t_0)$ versus $\widehat{H_{n,\alpha}^{IR2}}(t_0)$ for $t_0=0.3,0.5,0.7$ and $\widehat{H_{Y,2^{J_n},n}}(t_0)$ versus Standard Normal for $t_0=0.5$.}
\end{figure}
The QQ plots show that the probability distribution of our estimator for each $H(t_0)$ are close to that of IR2 estimator when the observed signal process is mBm and it is asymptotically normally distributed.

Through the above simulation study we conclude that there is no significant difference among IR2 estimator provided in Bardet and Surgailis (2013) and our wavelet-based estimator. And no significant difference is observed among wavelet-based estimators corresponding to different $\Phi\in C^2(\mathbb R)$. We also state that the bias and variance of $\widehat{H_{Y,2^{J_n},n}}(t)$ are generally greater than IR2 estimator when the sample size is relatively small. This is because, the wavelet-based method generally provides estimators of slower convergence rate than IR2 estimator.

\section{Appendix}
 \subsection{Proofs of (\ref{dnds}) and (\ref{dnd}).} First by using triangle inequality, we get
    \begin{eqnarray}
    \label{dnd1}
    &&|\widehat{d_{Y,n}}(2^{-j},k)-d_Y(2^{-j},k)|\nonumber\\
    &&\leq 2^{j/2}\sum_{l=0}^{2^{n-j}-1}\int_{\frac{l}{2^n}}^{\frac{l+1}{2^n}}|\psi(2^jt)||Y(l2^{-n}+k2^{-j})-Y(t+k2^{-j})|{\,\mathrm{d}} t.
    \end{eqnarray}
    Recall that $Y(t)=\Phi(X(t))$ for $t\in[0,1]$. Define the random variable $\|X\|_{\infty}$ to be
\begin{equation}
\label{supX}
\|X\|_{\infty}=\sup_{t\in[0,1]}|X(t)|.
\end{equation}
Since $\theta$ is continuous and not equal to $0$ almost everywhere, then $\{X(t)\}_{t\in[0,1]}$ is a Gaussian process with continuous trajectories, by applying Dudley's theorem and Borell's inequality (more precisely, with the same arguments for the proof of $\mathbb E(e^{\widetilde{V}})<+\infty$ on Page 1445-1446 in Rosenbaum (2008). See also Ledoux and Talagrand (2010)), we can show that
$
\mathbb E(e^{\|X\|_{\infty}})<+\infty.
$
This means all of $\|X\|_{\infty}$'s   moments are finite. Hence, using the mean value theorem, we get
\begin{equation}
\label{dnd2}
|Y(l2^{-n}+k2^{-j})-Y(t+k2^{-j})|\leq C_1|X(l2^{-n}+k2^{-j})-X(t+k2^{-j})|,
\end{equation}
where $C_1=\sup\limits_{s\in[-\|X\|_{\infty},\|X\|_{\infty}]}|\Phi'(s)|$ is a random variable. It follows from (\ref{dnd1}) and (\ref{dnd2}) that
 \begin{eqnarray}
    \label{dnd3}
    &&|\widehat{d_{Y,n}}(2^{-j},k)-d_Y(2^{-j},k)|\nonumber\\
    &&\leq C_12^{j/2}\sum_{l=0}^{2^{n-j}-1}\int_{\frac{l}{2^n}}^{\frac{l+1}{2^n}}|\psi(2^jt)||X(l2^{-n}+k2^{-j})-X(t+k2^{-j})|{\,\mathrm{d}} t.
    \end{eqnarray}
    In order to get (\ref{dnds}), we need (\ref{AH}), from which we see there exists a positive random variable $C_2$ with all finite moments such that
    \begin{eqnarray}
    \label{computeX}
    &&|X(l2^{-n}+k2^{-j})-X(t+k2^{-j})|\nonumber\\
    &&\le|\theta(l2^{-n}+k2^{-j})||B_{H(l2^{-n}+k2^{-j})}(l2^{-n}+k2^{-j})-B_{H(t+k2^{-j})}(t+k2^{-j})|\nonumber\\
    &&~~+|\theta(l2^{-n}+k2^{-j})-\theta(t+k2^{-j})||B_{H(t+k2^{-j})}(t+k2^{-j})|\nonumber\\
    &&\leq C_2 \big(|l2^{-n}-t|^{H(t+k2^{-j})}|\log|l2^{-n}-t||^{1/2}+|l2^{-n}-t|\big).
    \end{eqnarray}
     Observe that for $t\in[l2^{-n},(l+1)2^{-n}]$, $$|l2^{-n}-t|^{H(t+k2^{-j})}|\log|l2^{-n}-t||^{1/2}\ge |l2^{-n}-t|$$ and
     \begin{equation}
     \label{neighbor}\sup\limits_{{}^{n\ge0,j\leq n,l\leq 2^{n-j}-1,k\leq 2^j-1}_{l2^{-n}\le t\le (l+1)2^{-n}}}|l2^{-n}-t|^{H(t+k2^{-j})-H(k2^{-j})}|\log|l2^{-n}-t||^{1/2}<+\infty.
     \end{equation}
      This together with the fact that $|l2^{-n}-t|\leq 2^{-n}$ for $t\in[l2^{-n},(l+1)2^{-n}]$ yields there exists a positive random variable $C_3$ with all finite moments such that,
   \begin{equation}
   \label{dnds1}
    |X(l2^{-n}+k2^{-j})-X(t+k2^{-j})|\leq C_3 2^{-nH(2^{-j}k)}n^{1/2}.
    \end{equation}
    Then (\ref{dnds}) results from (\ref{dnd3}) and (\ref{dnds1}). $\square$

    Now we are going to prove (\ref{dnd}). For $r\ge1$, we consider the $r$-order moment of $|\widehat{d_{Y,n}}(2^{-j},k)-d_Y(2^{-j},k)|$ in (\ref{dnd3}). By applying the following two versions of Jensen's inequalities:
    \begin{equation}
    \label{Jensen}
    \Big(\sum_{i=1}^n|a_i|\Big)^r\leq n^{r-1}\sum_{i=1}^n|a_i|^r~\mbox{and}~\Big(\int_a^b|f(s)|{\,\mathrm{d}} s\Big)^r\leq |b-a|^{r-1}\int_a^b|f(s)|^r{\,\mathrm{d}} s,
   \end{equation}
     and Cauchy-Schwarz inequality, we obtain
     \begin{eqnarray}
    \label{dnd4}
    &&\mathbb E|\widehat{d_{Y,n}}(2^{-j},k)-d_Y(2^{-j},k)|^r\nonumber\\
    &&\leq 2^{jr/2}2^{-j(r-1)}\!\!\!\!\sum_{l=0}^{2^{n-j}-1}\!\!\!\!\int_{\frac{l}{2^n}}^{\frac{l+1}{2^n}}|\psi(2^jt)|^r\mathbb E\big(C_1|X(l2^{-n}+k2^{-j})-X(t+k2^{-j})|\big)^r{\,\mathrm{d}} t\nonumber\\
    &&\leq \big(\sup_{s\in[0,1]}|\psi(s)|^r\big)2^{-j(r/2-1)}\nonumber\\
    &&~~\times\sum_{l=0}^{2^{n-j}-1}\int_{\frac{l}{2^n}}^{\frac{l+1}{2^n}}\big(\mathbb E(C_1^{2r})\big)^{1/2}\big(\mathbb E|X(l2^{-n}+k2^{-j})-X(t+k2^{-j})|^{2r}\big)^{1/2}{\,\mathrm{d}} t.
    \end{eqnarray}
    Note that by Lemma 2.12 (i) in Ayache et al. (2011), there exists a constant $c_1>0$ which does not depend on $n,l,j,k$ and $H$ such that
    \begin{eqnarray*}
    &&\mathbb E|B_{H(l2^{-n}+k2^{-j})}(l2^{-n}+k2^{-j})-B_{H(t+k2^{-j})}(t+k2^{-j})|^2\\
    &&\leq c_1|l2^{-n}+k2^{-j}-(t+k2^{-j})|^{2\max\{H(l2^{-n}+k2^{-j}),H(t+k2^{-j})\}}\\
    &&\leq c_1|l2^{-n}+k2^{-j}-(t+k2^{-j})|^{2H(l2^{-n}+k2^{-j})}\\
    &&=c_1|l2^{-n}-t|^{2H(k2^{-j})}|l2^{-n}-t|^{2H(l2^{-n}+k2^{-j})-2H(k2^{-j})}.
    \end{eqnarray*}
    Therefore by using again (\ref{neighbor}) and the fact that $|l2^{-n}-t|\leq 2^{-n}$, there exists some constant $c>0$ such that
    \begin{equation}
    \label{dnd5'}
    \mathbb E|B_{H(l2^{-n}+k2^{-j})}(l2^{-n}+k2^{-j})-B_{H(t+k2^{-j})}(t+k2^{-j})|^2\leq c2^{-2nH(k2^{-j})}.
    \end{equation}
    Using (\ref{dnd5'}) and similar computations as in (\ref{computeX}), we obtain there exists $c_2>0$ such that
    \begin{equation}
    \label{dnd5}
    \mathbb E|X(l2^{-n}+k2^{-j})-X(t+k2^{-j})|^2\leq c_22^{-2nH(k2^{-j})}.
    \end{equation}
    By using the fact that all the moments of Gaussian variable are equivalent, we get there exists some constant $c_3>0$ (only depending on $r$) such that
    \begin{equation}
    \label{dnd6}
    \mathbb E|X(l2^{-n}+k2^{-j})-X(t+k2^{-j})|^{2r}\leq c_32^{-2rnH(k2^{-j})}.
    \end{equation}
    Finally it results from (\ref{dnd4}) and (\ref{dnd6}) that
    \begin{eqnarray*}
    \mathbb E|\widehat{d_{Y,n}}(2^{-j},k)-d_Y(2^{-j},k)|^r\leq c 2^{-r(nH(2^{-j}k)+j/2)},
    \end{eqnarray*}
    where $c=\big(\sup_{s\in[0,1]}|\psi(s)|^r\big)\big(c_3\mathbb E(C_1^{2r})\big)^{1/2}$. Therefore (\ref{dnd}) has been proven. $\square$
    \subsection{Proof of (\ref{wideV})}
     First notice that, by the definition of $\nu_{t_0,2^j}$, we have
    \begin{equation}
    \label{card}
    |card(\nu_{t_0,2^j})- 2^{j+1}\epsilon_j|\leq 3
    \end{equation}
    as $j\rightarrow+\infty$, because $2^{j}\epsilon_j\ge1$.

    It follows from (\ref{Jensen}), Cauchy-Schwarz inequality, (\ref{card}) and the fact that $(a+b)^4\le 2^3(a^4+b^4)$ that
    \begin{eqnarray}
    \label{VnV1}
    && \mathbb E|\widehat{V_{n,t_0,j}}-V_{Y,t_0,j}|^2\leq card(\nu_{t_0,2^j})\sum_{k\in\nu_{t_0,2^j}}\mathbb E\big|\widehat{d_{Y,n}}(2^{-j},k)^2-d_Y(2^{-j},k)^2\big|^2\nonumber\\
    &&\leq 3\times2^j\epsilon_j\sum_{k\in\nu_{t_0,2^j}}2^{3/2}(\mathbb E|\widehat{d_{Y,n}}(2^{-j},k)|^4+\mathbb E|d_Y(2^{-j},k)|^4)^{1/2}\nonumber\\
    &&~~~~\times(\mathbb E|\widehat{d_{Y,n}}(2^{-j},k)-d_Y(2^{-j},k)|^4)^{1/2}.
    \end{eqnarray}
    Roughly speaking (and it can be proven without efforts), since the trajectory $\{\Phi(X(t))\}_{t\ge0}$ is at least as smooth as $\{X(t)\}_{t\ge0}$, then for $r\ge1$, there exists a constant $c_4>0$ (only depending on $r$) such that,
    \begin{eqnarray}
    \label{VnV2}
    &&\mathbb E|\widehat{d_{Y,n}}(2^{-j},k)|^r\leq c_42^{-jr(H(k2^{-j})+1/2)};\nonumber\\
    &&\mathbb E|d_Y(2^{-j},k)|^r\leq c_42^{-jr(H(k2^{-j})+1/2)}.
    \end{eqnarray}
    Then it results from (\ref{VnV1}), (\ref{VnV2}), (\ref{dnd}) and (\ref{card}) that
    \begin{equation}
    \label{VnV3}
    \mathbb E|\widehat{V_{n,t_0,j}}-V_{Y,t_0,j}|^2\leq 9\times2^{2j}\epsilon_j^2\big((8c_4)^{1/2}2^{-j(2H(k2^{-j})+1)}\big)\big(c^{1/2}2^{-2nH(k2^{-j})-j}\big).
    \end{equation}
  In view of the equivalence relation between $H(k2^{-j})$ and $H(t_0)$ as $k\in\nu_{t_0,2^j}$ and $j\rightarrow+\infty$,
    (\ref{wideV}) finally results from (\ref{VnV3}). $\square$
\subsection{Proof of Proposition \ref{prop1}}
\subsubsection{Proof of (\ref{prop11})}

By using the fact that $d_X(a,k)$, $d_X(b,k')$ are zero-mean Gaussian random variables, Fubini's theorem, the isometry property of mBm's harmonizable presentation and a change of variables, we get
\begin{eqnarray}
\label{cov1}
&&Cov(d_X( a,k) ,d_X(b,k'))=\frac{1}{\sqrt{ab}}\int_{0}^{a}\int_{0}^{b}\psi \big( \frac{t}{a}%
\big) \psi \big( \frac{s}{b}\big) \mathbb E\big( X(t+ak)X(s+bk^{\prime })%
\big) {\,\mathrm{d}} t{\,\mathrm{d}} s\nonumber\\
&&=\sqrt{ab}\int_{0}^{1}\int_{0}^{1}\psi ( t%
) \psi ( s)\theta(at+ak)\theta(bs+bk')\nonumber\\
&&~~~~\times \int_{\mathbb{R}}\frac{\big(
e^{i(at+ak) u}-1\big) \big( e^{-i( bs+bk')
u}-1\big) }{|u| ^{H( at+ak) +H(
bs+bk') +1}}{\,\mathrm{d}} u{\,\mathrm{d}} t{\,\mathrm{d}} s.
\end{eqnarray}
Since the pointwise H\"{o}lder exponent of $B_{H(t)}(t)$ in the neighborhood of $ak,bk'$ behave
locally asymptotically like those of fractional Brownian motions with Hurst parameters $H(ak),H(bk')$: $B_{H(ak)}(t)$ on $ak$ and $%
B_{H( bk')}(t)$ on $bk'$, we can thus consider a Taylor expansion
of $H$ respectively on $ak$
and $bk^{\prime }$. To be more explicit, let's fix $u\in\mathbb R$ and define $f(x) =\theta(x)|u| ^{-H(x) -%
1/2}$, since $f$ belongs to $C^2([0,1])$, we take the second order Taylor expansion for $f$ respectively on $%
ak$ and $bk^{\prime }.$ There exist $\xi_t \in ( ak,ak+at) $ and $%
\xi ^{\prime }_s\in ( bk^{\prime },bk^{\prime }+bs) $ such that
\begin{eqnarray}
\label{A}
f( at+ak) &=&A_{0}(u,ak) +A_{1}(u,at,ak)+A_2(u,at,\xi_t); \\
\label{A'}
f( bs+bk')
&=&A_{0}( u,bk') +A_{1}( u,bs,bk')+A_2(u,bs,\xi'_s),
\end{eqnarray}%
where we denote, for $u\neq 0$, $x,y\ge0$,
\begin{eqnarray*}
&&A_0(u,y)=\theta(y)|u| ^{-H(y) -1/2};\\
&&A_{1}( u,x,y ) =x\big(\theta'(y)|u|^{-H(y)-1/2}-\theta(y)H'(y)|u|^{-H(y)-1/2}\log|u|\big);
\end{eqnarray*}
and
\begin{eqnarray*}
&&A_2(u,x,y)=\frac{1}{2}x^2\big(\theta''(y)|u|^{-H(y)-1/2}-(2\theta'(y)H'(y)+\theta(y)H''(y))\\
&&~~\times|u|^{-H(y)-1/2}\log|u|+\theta(y)\big(H'(y)\big)^2|u|^{-H(y)-1/2}\big(\log|u|\big)^2\big).
\end{eqnarray*}
Thus we rewrite (\ref{cov1}) as
$$Cov\big(d_X( a,k) ,d_X(b,k')\big)=\sum\limits_{l,l^{\prime }\in \left\{ 0,1,2\right\} }\mathcal{I}%
_{l,l^{\prime }}\left( k,k^{\prime },a,b\right),
$$
where
\begin{eqnarray*}
&&\mathcal{I}_{l,l^{\prime }}( k,k^{\prime },a,b)\\
&& =\sqrt{ab}%
\int_{0}^{1}\int_{0}^{1}\int_{\mathbb{R}}\psi (t) \psi
(s) \big( e^{i\left( at+ak\right) u}-1\big) \big(
e^{-i\left( bs+bk^{\prime }\right) u}-1\big) A_{l}A_{l^{\prime }}{\,\mathrm{d}} u{\,\mathrm{d}} t{\,\mathrm{d}} s.
\end{eqnarray*}%
We note here $A_{l}$'s and $A_l'$'s are notations in short for (\ref{A}) and (\ref{A'}).
By using (\ref{cov1}), it suffices to make an identification of all the terms $\mathcal{I}_{l,l^{\prime
}}\left( k,k^{\prime },a,b\right) ^{\prime }s$ in order to estimate the covariance structure of the wavelet coefficients. We consider different cases according to the values of $(l,l')$. The key to these identifications is to observe the following:\begin{itemize}
                        \item First, observe that for $x,y>0,\ x\neq y,\ \alpha >0$, $p\in\mathbb N$, we have
\begin{eqnarray}
\label{general}
&&\int_{\mathbb{R}}\frac{( e^{ixu}-1)(
e^{-iyu}-1) }{\left\vert u\right\vert ^{\alpha+1 }}(\log|u|)^p {\,\mathrm{d}} u\nonumber\\
&&=\frac{1}{2}\sum_{l=0}^p(-1)^lC_{p-l}(\alpha)\big( \left\vert x\right\vert ^{\alpha}(\log \left\vert
x\right\vert)^l +\left\vert y\right\vert ^{\alpha}(\log \left\vert
y\right\vert)^l -\left\vert x-y\right\vert ^{\alpha }(\log \left\vert
x-y\right\vert)^l \big),\nonumber\\
\end{eqnarray}
where  for $l\in\{0,\ldots,p\}$, $C_{l}( \alpha ) =\big({}_l^p\big)\int_{\mathbb{R}}\frac{| e^{iu}-1| ^{2}}{%
\left\vert u\right\vert ^{\alpha+1 }} (\log
\left\vert u\right\vert)^l {\,\mathrm{d}} u$, with $\big({}_l^p\big)=\frac{p!}{l!(p-l)!}$ being the binomial coefficient.
                        \item Secondly, for $|\frac{at-bs}{ak-bk'}|\leq1$, $Q\ge1$, $l,l'\in\{0,\ldots,Q\}$, $p\in\mathbb N$ and $\alpha>0$, a $2Q-l-l'$ order Taylor expansion of $q_{\alpha,p}(x)=(1+x)^{\alpha}(\log|1+x|)^p$ on $x=\frac{at-bs}{ak-bk'}$ yields:
                        \begin{eqnarray}
\label{computeI}
&&%
\int_{0}^{1}\int_{0}^{1}t^ls^{l'}\psi (t) \psi (s)
|at+ak-bs-bk'| ^{\alpha }(\log|at+ak-bs-bk'|)^p{\,\mathrm{d}} t{\,\mathrm{d}} s\nonumber\\
&&=|ak-bk'| ^{\alpha }\sum_{j=0}^p\big({}_j^p\big)(\log|ak-bk'|)^{p-j}\int_{0}^{1}\!\!\!\!\int_{0}^{1}t^ls^{l'}\psi
(t) \psi (s) \Big( 1+\frac{at-bs}{ak-bk'}\Big) ^{\alpha }\nonumber\\
&&~~~~\times\Big(\log\big| 1+\frac{at-bs}{ak-bk'}\big|\Big)^{j}{\,\mathrm{d}} t{\,\mathrm{d}} s\nonumber\\
&&=|ak-bk'| ^{\alpha}\sum_{j=0}^p\big({}_j^p\big)(\log|ak-bk'|)^{p-j}\int_{0}^{1}\int_{0}^{1}t^ls^{l'}\psi
(t) \psi (s)\nonumber\\
&&~~~~\times \theta_{\alpha,2Q-l-l',j}\big(\frac{at-bs}{ak-bk'}\big)\Big(\frac{at-bs}{ak-bk'}\Big) ^{2Q-l-l'}{\,\mathrm{d}} t{\,\mathrm{d}} s\nonumber\\
&&=\frac{(ab) ^{Q-(l+l')/2}}{%
\left\vert ak-bk^{\prime }\right\vert ^{2Q-l-l'-\alpha}}\Big(\sum_{j=0}^p\big({}_j^p\big)A_{\alpha,2Q-l-l',j}\big(\frac{a}{b}\big)(\log|ak-bk'|)^{p-j}\Big),
\end{eqnarray}
where the integral remainder $\theta_{\alpha,2Q-l-l',j}$ of $q(\cdot)$'s $(2Q-l-l')$-th order Taylor expansion (see e.g. Apostol 1967) is given as: for $2Q-l-l'=0$,
\begin{equation}
\label{theta'}
\theta_{\alpha,2Q-l-l',j}\big(\frac{at-bs}{ak-bk'}\big)=q_{\alpha,j}\big(\frac{at-bs}{ak-bk'}\big);
\end{equation}
for $2Q-l-l'\ge1$,
\begin{eqnarray}
\label{theta}&&\theta_{\alpha,2Q-l-l',j}\Big(\frac{at-bs}{ak-bk'}\Big)\nonumber\\
&& = \frac{1 }{%
(2Q-l-l'-1) !}\int_0^1\!\!(1-\eta)^{2Q-l-l'-1}q_{\alpha,j}^{(2Q-l-l')}\Big(\eta \frac{at-bs}{ak-bk^{\prime }}+(1-\eta)\Big){\,\mathrm{d}}\eta;\nonumber\\
\end{eqnarray}
and the function
\begin{eqnarray}
\label{Aab}
&&A_{\alpha,2Q-l-l',j}\big(\frac{a}{b}\big)\nonumber\\
&&=\int_{0}^{1}\int_{0}^{1}t^ls^{l'}\psi (t)
\psi (s)\theta_{\alpha,2Q-l-l',j}\big(\frac{at-bs}{ak-bk'}\big)\big(\sqrt{\frac{a}{b}}t-\sqrt{\frac{b}{a}}s\big)^{2Q-l-l'}{\,\mathrm{d}} t{\,\mathrm{d}} s.\nonumber\\
\end{eqnarray}
Here we note that for $r\in\mathbb N$ and $|x|< 1$,
$$
q_{\alpha,p}^{(r)}(x)=(1+x)^{\alpha-r}\!\!\!\!\sum_{i,j\in\mathbb N,i+j=r}\!\!\!\!(\alpha-1)\ldots(\alpha-i)p\ldots(p-j)(\log|1+x|)^{p-j}.
$$
                      \end{itemize}
\begin{description}
\item[Case (i)] $l=l^{\prime }=0 $.

In this case we have
\begin{eqnarray*}
&&\mathcal{I}_{0,0}=\left( k,k^{\prime
},a,b\right)\theta(ak)\theta(bk')\sqrt{ab}\\
&&~~~~\times\int_{0}^{1}\int_{0}^{1}\int_{\mathbb{R}}\psi \left( t\right)
\psi \left( s\right) \frac{( e^{i\left( at+ak\right) u}-1)
( e^{-i\left( bs+bk^{\prime }\right) u}-1) }{\left\vert
u\right\vert ^{H\left( ak\right) +H\left( bk^{\prime }\right) +1}}{\,\mathrm{d}} u{\,\mathrm{d}} t{\,\mathrm{d}} s.
\end{eqnarray*}
Let $p=0$, $x=at+ak,\ y=bs+bk^{\prime }$ and $\alpha =H\left( ak \right)
+H\left( bk^{\prime }\right) $ in (\ref{general}). It follows
\begin{eqnarray*}
&&\mathcal{I}_{0,0}\left( k,k^{\prime },a,b\right)=-\frac{C_0(H(ak)+H(bk'))}{2}\theta(ak)\theta(bk')\sqrt{ab}\\
&&~~~~\times\int_{0}^{1}\int_{0}^{1}\psi (t) \psi (s) \left\vert at+ak-bs-bk^{\prime }\right\vert ^{H\left( ak\right) +H\left(
bk^{\prime }\right) }{\,\mathrm{d}} t{\,\mathrm{d}} s.
\end{eqnarray*}%
Then by the assumption $\sup\limits_{t,s\in[0,1]}\vert \frac{at-bs}{ak-bk^{\prime }}\vert \leq 1$, we can thus take $l=l'=0$ in (\ref{computeI}) and obtain
\begin{eqnarray*}
&&\mathcal{I}_{0,0}\left( k,k^{\prime },a,b\right)\nonumber\\ &&=-\frac{C_0(H(ak)+H(bk'))}{2}A_{H(ak)+H(bk'),2Q,0}\big(\frac{a}{b}\big)\frac{(ab) ^{Q+1/2}\theta(ak)\theta(bk')}{%
\left\vert ak-bk^{\prime }\right\vert ^{2Q-H(ak)-H(bk')}}.\nonumber\\
\end{eqnarray*}
We finally obtain
\begin{equation}
\label{computeI00}
\mathcal{I}_{0,0}\left( k,k^{\prime },a,b\right) =C_1(H(ak)+H(bk'),Q,\frac{a}{b})\frac{\left( ab\right)
^{Q+1/2}\theta(ak)\theta(bk')}{\left\vert ak-bk^{\prime }\right\vert ^{2Q-H(ak)-H(bk') }},
\end{equation}%
where the function $$C_1(H(ak)+H(bk'),Q,\frac{a}{b})=-\frac{C_0(H(ak)+H(bk'))}{2}A_{H(ak)+H(bk'),2Q,0}.$$
\item[Case (ii)] $(l,l')\in\{(1,0),(0,1)\}.$

By definition, $\mathcal{I}_{1,0}( k,k^{\prime },a,b)$
equals
\begin{eqnarray*}
&&\theta'(ak)\theta(bk')a^{\frac{3}{2}}b^{\frac{1}{2}}\int_{0}^{1}\!\!\!\int_{0}^{1}\!\!\!\int_{\mathbb{R}} t\psi(t) \psi (s) \frac{(
e^{i( at+ak) u}-1) ( e^{-i\left( bs+bk^{\prime
}\right) u}-1 }{|u| ^{H(ak)
+H(bk') +1}} {\,\mathrm{d}} u{\,\mathrm{d}} t{\,\mathrm{d}} s\\
&&~~-\theta(ak)\theta(bk')H'(ak)a^{\frac{3}{2}}b^{\frac{1}{2}}\\
&&~~~~\times\int_{0}^{1}\!\!\!\int_{0}^{1}\!\!\!\int_{\mathbb{R}} t\psi(t) \psi (s) \frac{(
e^{i( at+ak) u}-1) ( e^{-i\left( bs+bk^{\prime
}\right) u}-1) }{|u| ^{H(ak)
+H(bk') +1}}\log |u| {\,\mathrm{d}} u{\,\mathrm{d}} t{\,\mathrm{d}} s.
\end{eqnarray*}
Since $\theta,H\in C^{2}([0,1])$, then by the triangle inequality the above item can be expressed as
\begin{eqnarray*}
&&\sum_{j=0}^1\mathcal O\Big(a^{\frac{3}{2}}b^{\frac{1}{2}}\\
&&~~\times\int_{0}^{1}\!\!\!\int_{0}^{1}\!\!\!\int_{\mathbb{R}}t\psi(t) \psi (s) \frac{(
e^{i( at+ak) u}-1)( e^{-i( bs+bk^{\prime
}) u}-1) }{|u| ^{H(ak)
+H(bk') +1}}(\log |u|)^j {\,\mathrm{d}} u{\,\mathrm{d}} t{\,\mathrm{d}} s\Big).
\end{eqnarray*}
 Let $p=0$ and $p=1$ (respectively corresponding to $j=0$ and $1$ of the above expression), $x=at+ak,\ y=bs+bk'$ and $\alpha =H(ak)
+H( bk')$ in (\ref{general}), $(l,l')=(1,0)$ in (\ref{computeI}), we obtain
\begin{equation*}
\label{computeI10}
\mathcal{I}_{1,0}\left( k,k^{\prime },a,b\right)=\mathcal O\Big(\frac{(ab)^{Q}a}{|ak-bk'|^{2Q-H(ak)-H(bk')-1}}\big(1+|\log|ak-bk'||\big)\Big).
\end{equation*}
Since  $\mathcal{I}_{1,0}( k,k^{\prime },a,b)=\mathcal{I}_{0,1}( k',k,b,a)$, thus we get at the meanwhile,
\begin{equation*}
\mathcal{I}_{0,1}\left( k,k^{\prime },a,b\right)=\mathcal O\Big(\frac{(ab)^{Q}b}{|ak-bk'|^{2Q-H(ak)-H(bk')-1}}\big(1+|\log|ak-bk'||\big)\Big).
\end{equation*}
Using the fact that $a\sim b$, we conclude
\begin{eqnarray}
\label{computeI01}
&&\mathcal{I}_{1,0}\left( k,k',a,b\right)+\mathcal{I}_{0,1}\left( k,k^{\prime },a,b\right)\nonumber\\
&&=\mathcal O\Big(\frac{(ab)^{Q+1/2}}{|ak-bk'|^{2Q-H(ak)-H(bk')-1}}\big(1+|\log|ak-bk'||\big)\Big).
\end{eqnarray}
\item[Case (iii)] $l=l'=1.$

Observe that $\mathcal{I}_{1,1}\left( k,k^{\prime },a,b\right)$ can be expressed as
\begin{eqnarray*}
&&\sum_{j=0}^2\mathcal O\Big(a^{\frac{3}{2}}b^{\frac{3}{2}}\\
&&~~\times\int_{0}^{1}\!\!\!\int_{0}^{1}\!\!\!\int_{\mathbb{R}}st\psi(t) \psi (s) \frac{\big(
e^{i( at+ak) u}-1\big) \big( e^{-i\left( bs+bk^{\prime
}\right) u}-1\big) }{|u| ^{H(ak)
+H(bk') +1}}(\log |u|)^j {\,\mathrm{d}} u{\,\mathrm{d}} t{\,\mathrm{d}} s\Big).
\end{eqnarray*}

 Let $p=0,1,2$ respectively, $x=at+ak,\ y=bs+bk^{\prime }$ and $\alpha =H\left( ak \right)
+H\left( bk^{\prime }\right)$ in (\ref{general}) and $(l,l')=(1,1)$ in (\ref{computeI}), we get
\begin{equation}
\label{computeI11}
\mathcal{I}_{1,1}\left( k,k^{\prime },a,b\right)=\sum_{j=0}^2\mathcal O\Big(\frac{(ab)^{Q+1/2}}{|ak-bk'|^{2Q-H(ak)-H(bk')-2}}|\log|ak-bk'||^j\Big).
\end{equation}
\item[Case (iv)] $(l,l')\in\{(2,0),(0,2)\}.$\\
By using the fact that $\theta,H\in C^2([0,1])$, $\mathcal{I}_{2,0}\left( k,k^{\prime },a,b\right)$ can be expressed as
$$
\sum_{p=0}^2\!\mathcal O\Big(a^{\frac{5}{2}}b^{\frac{1}{2}}\!\!\int_{0}^{1}\!\!\Big|\!\int_{0}^{1}\!\!\!\!\int_{\mathbb{R}} \psi (s) \frac{(
e^{i( at+ak) u}-1) ( e^{-i( bs+bk^{\prime
}) u}-1) }{|u| ^{H(\xi_t)
+H(bk') +1}}(\log|u|)^p {\,\mathrm{d}} u{\,\mathrm{d}} s\Big|{\,\mathrm{d}} t\Big).
$$
 Let $p=0,1,2$ respectively, $x=at+ak,\ y=bs+bk^{\prime }$ and $\alpha =H\left( \xi_t \right)
+H\left( bk^{\prime }\right)$ in (\ref{general}) and notice that $\xi_t$ depends on $t$, we get
\begin{eqnarray*}
&&\mathcal{I}_{2,0}\left( k,k^{\prime },a,b\right)=\sum_{p=0}^2\mathcal O\Big(a^{\frac{5}{2}}b^{\frac{1}{2}}\int_{0}^{1}\Big|\int_{0}^{1}\int_{\mathbb{R}} \psi (s)\sum_{l=0}^p(-1)^lC_{p-l}(\alpha)\\
&&~~~~\times\Big( \left\vert bs+bk'\right\vert ^{\alpha}(\log \left\vert
bs+bk'\right\vert)^l\\
&&~~~~-|at-bs+ak-bk'| ^{\alpha }(\log |at-bs+ak-bk'| )^l\Big) {\,\mathrm{d}} u{\,\mathrm{d}} s\Big|{\,\mathrm{d}} t\Big).
\end{eqnarray*}

Then similarly to (\ref{computeI}), using a $Q$ order Taylor expansion of $q_{\alpha,p}(\cdot)$ respectively on $s/k'$ and on $(at-bs)/(ak-bk')$, and also use the fact that for $x>0$, $x^{H(\xi_t)}\sim x^{H(ak)}$ (since $\xi_t\in(ak,ak+at)$), we obtain

\begin{eqnarray}
\label{computeI20}
&&\mathcal{I}_{2,0}\left( k,k^{\prime },a,b\right)=\mathcal O\Big(\frac{a^{\frac{5}{2}}b^{\frac{1}{2}+H(ak)+H(bk')}}{|k'|^{Q-H(ak)-H(bk')}}\big(1+|\log|bk'||+|\log|bk'||^2\big)\Big)\nonumber\\
&&~~+\mathcal O\Big(\frac{a^{\frac{5}{2}}b^{\frac{1}{2}}(ab)^{\frac{Q}{2}}}{|ak-bk'|^{Q-H(ak)-H(bk')}}\big(1+|\log|ak-bk'||+|\log|ak-bk'||^2\big)\Big).\nonumber\\
\end{eqnarray}
By using symmetric property,
\begin{eqnarray}
\label{computeI02}
&&\mathcal{I}_{0,2}\left( k,k^{\prime },a,b\right)=\mathcal O\Big(\frac{b^{\frac{5}{2}}a^{\frac{1}{2}+H(ak)+H(bk')}}{|k|^{Q-H(ak)-H(bk')}}\big(1+|\log|ak||+|\log|ak||^2\big)\Big)\nonumber\\
&&~~+\mathcal O\Big(\frac{b^{\frac{5}{2}}a^{\frac{1}{2}}(ab)^{\frac{Q}{2}}}{|ak-bk'|^{Q-H(ak)-H(bk')}}\big(1+|\log|ak-bk'||+|\log|ak-bk'||^2\big)\Big).\nonumber\\
\end{eqnarray}
\item[Case (v)] $l+l'\ge3$.\\
Since the computations are quite similar as in the previous cases, we present the results without proof.
\begin{eqnarray}
\label{computeI21}
&&\mathcal{I}_{2,1}\left( k,k^{\prime },a,b\right)=\mathcal O\Big(\frac{a^{\frac{5}{2}}b^{\frac{3}{2}+H(ak)+H(bk')}}{|k'|^{Q-H(ak)-H(bk')-1}}\sum_{j=0}^3|\log|bk'||^j\nonumber\\
&&~~~~+\frac{a^{2}b(ab)^{\frac{Q}{2}}}{|ak-bk'|^{Q-H(ak)-H(bk')-1}}\sum_{j=0}^3|\log|ak-bk'||^j\Big).
\end{eqnarray}
\begin{eqnarray}
\label{computeI12}
&&\mathcal{I}_{1,2}\left( k,k^{\prime },a,b\right)=\mathcal O\Big(\frac{b^{\frac{5}{2}}a^{\frac{3}{2}+H(ak)+H(bk')}}{|k|^{Q-H(ak)-H(bk')-1}}\sum_{j=0}^3|\log|ak||^j\nonumber\\
&&~~~~+\frac{b^{2}a(ab)^{\frac{Q}{2}}}{|ak-bk'|^{Q-H(ak)-H(bk')-1}}\sum_{j=0}^3|\log|ak-bk'||^j\Big).
\end{eqnarray}
\begin{eqnarray}
\label{computeI22}
&&\mathcal{I}_{2,2}\left( k,k^{\prime },a,b\right)=\mathcal O\Big(\frac{a^{\frac{5}{2}}b^{\frac{5}{2}+H(ak)+H(bk')}}{|k'|^{-H(ak)-H(bk')}}\sum_{j=0}^4|\log|bk'||^j\nonumber\\
&&~~~~+\frac{b^{\frac{5}{2}}a^{\frac{5}{2}+H(ak)+H(bk')}}{|k|^{-H(ak)-H(bk')}}\sum_{j=0}^4|\log|ak||^j\nonumber\\
&&~~~~+\frac{(ab)^{\frac{5}{2}}}{|ak-bk'|^{-H(ak)-H(bk')}}\sum_{j=0}^4|\log|ak-bk'||^j\Big).
\end{eqnarray}
\end{description}
Now denote by
\begin{equation}
\label{Tkk'}
T(k,k',Q,a,b)=(ab)^{-\frac{H(ak)+H(bk')}{2}-\frac{1}{2}}\sum_{l,l'\in\{0,1,2\},(l,l')\neq(0,0)}\mathcal{I}_{l,l'}(k,k',a,b).
\end{equation}
It remains to show that for $Q\ge2$ and $a\sim b$,
$$
\sum_{k=0}^{[a^{-1}]-1}\sum_{k'=0}^{[b^{-1}]-1}\big(T(k,k',Q,a,b)\big)^2=\mathcal O((ab)^{1/2}\log a\log b),~\mbox{as $a,b\rightarrow0$}.
$$
According to (\ref{Tkk'}), it suffices to prove for any $(l,l')\neq (0,0)$,
\begin{equation}
\label{boundI}
\sum_{k=0}^{[a^{-1}]-1}\sum_{k'=0}^{[b^{-1}]-1}(ab)^{-H(ak)-H(bk')-1}\big(\mathcal{I}_{l,l'}(k,k',a,b)\big)^2=\mathcal O((ab)^{1/2}\log a\log b).
\end{equation}
To  prove (\ref{boundI}) holds we only take $\mathcal I_{1,0}(k,k',a,b)$ as an example, since the computations of the other items are similar. Recall that
$$
\frac{\big(\mathcal I_{1,0}(k,k',a,b)\big)^2}{(ab)^{H(ak)+H(bk')+1}}=\mathcal O\Big(\frac{ab\big(1+|\log (ab)|+|\log|k\sqrt{a/b}-k'\sqrt{b/a}||\big)^2}{|k\sqrt{a/b}-k'\sqrt{b/a}|^{4Q-2H(ak)-2H(bk')-2}}\Big).
$$
    Remember that, the fact that $Q\ge2$ yields $4Q-2H(ak)-2H(bk')-2>2$, and $|k\sqrt{a/b}-k'\sqrt{b/a}|\ge\big[|k\sqrt{a/b}-k'\sqrt{b/a}|\big]$. These facts together with $a\sim b$ imply
\begin{eqnarray*}
&&\sum_{k=0}^{[a^{-1}]-1}\sum_{k'=0}^{[b^{-1}]-1}(ab)^{-H(ak)-H(bk')-1}\big(\mathcal I_{1,0}(k,k',a,b)\big)^2\\
&&=\mathcal O\Big(ab|\log(ab)|^2\sum_{k=0}^{[a^{-1}]-1}\sum_{|l|=1}^{+\infty}\frac{(\log|l|)^2}{|l|^{4Q-4\sup_{t\in[0,1]}H(t)-2}}\Big)\\
&&=\mathcal O\big(a^{-1}ab|\log(ab)|^2\big)=\mathcal O\big(b|\log(ab)|^2\big)=\mathcal O\big((ab)^{1/2}\log a\log b\big).~~\square\\
\end{eqnarray*}
\subsubsection{Proof of (\ref{prop12})}
When $ak=bk'$, we still make a Taylor expansion as in the proof of (\ref{prop11}).
For the first term $\mathcal{I}_{0,0}( k,k',a,b) $, since $b/a=k/k'=\varrho$ (some constant), then using the same formula in (\ref{computeI00}), we get%
\begin{eqnarray*}
\mathcal{I}_{0,0}\big( k,k',a,b\big)&=&-a\int_{0}^{1}\int_{0}^{1}\psi
\left( t\right) \psi \left( s\right) \frac{C_{0}(2H(ak))}{2}\theta(ak)^2\left\vert
at-\varrho as\right\vert ^{2H\left( ak\right) }{\,\mathrm{d}} t{\,\mathrm{d}} s \\
&=&C_2(ak,\varrho)a^{2H\left( ak\right) +1},
\end{eqnarray*}
where
\begin{equation}
\label{c2}
C_2(ak,\varrho)=-\theta(ak)^2\frac{C_{0}(2H(ak))}{2}\int_{0}^{1}\int_{0}^{1}\psi
\left( t\right) \psi \left( s\right) \left\vert t-\varrho s\right\vert ^{2H\left(
ak\right) }{\,\mathrm{d}} t{\,\mathrm{d}} s.
\end{equation}

The remaining items of $\mathcal I_{l,l'}$'s are of higher order so we only give a bound of them. By means of a Taylor expansion of $r_{\alpha,p}(x)=x^{\alpha}(\log|x|)^p$ on $t-\varrho s$, then similar discussion shows that, for $l+l'\ge1$,
$$
\mathcal{I}_{l,l'}\big( k,k',a,b\big)=\mathcal{I}_{l,l'}\big( k,k/\varrho,a,\varrho a\big)=\mathcal O\big(a^{2H(
ak) +2}|\log a|^4\big).~~\square
$$
\subsection{Proof of Proposition \ref{VarVX}}
Fix $t_0\in(0,1)$. It follows from the definition of $V_{X,t_0,j}$, (\ref{VXV}), (\ref{lemme12}) and the fact that $\limsup\limits_{j\rightarrow+\infty,k\in\nu_{t_0,2^j}}k2^{-j}=t_0$ that
\begin{eqnarray}
\label{VarV'1}
&&\mathbb E\big(V_{X,t_0,j}\big)=\sum_{k\in\nu_{t_0,2^j}}c_2(2^{-j}k)2^{-j(2H(2^{-j}k) +1)}+\mathcal{O%
}\Big( \sum_{k\in\nu_{t_0,2^j}}2^{-2j(H(2^{-j}k) +1)}j^4\Big)\nonumber\\
&&=\sum_{k\in\nu_{t_0,2^j}}c_2(2^{-j}k)2^{-j(2H(2^{-j}k) +1)}+\mathcal{O%
}\big( 2^{-j(2H(t_0) +1)}j^4\epsilon_j\big).
\end{eqnarray}
It remains to show, for $k\in\nu_{t_0,2^j}$,
$$
c_2(2^{-j}k)2^{-j(2H(2^{-j}k) +1)}=c_2(t_0)2^{-j(2H(t_0) +1)}+\mathcal O(2^{-j(2H(t_0) +1)}j\epsilon_j).
$$
For this purpose, for fixed $j,k$ we define $K_1\in C^1([0,1])$ as for $x\in[0,1]$,
$$K_1(x)=c_2(x)2^{-j(2H(x) +1)}.
$$
Then by using the Taylor expansion of $K_1$ around $x=t_0$ and the facts that $|2^{-j}k-t_0|=\mathcal O(\epsilon_j)$, $c_2\in C^2([0,1])$, $H\in C^2([0,1])$, we obtain
\begin{eqnarray}
\label{VarV'2}
&&K_1(2^{-j}k)=K_1(t_0)+K_1'(\xi)(2^{-j}k-t_0)=K_1(t_0)\nonumber\\
&&~~+\big(c_2'(\xi)2^{-j(2H(\xi) +1)}+c_2(\xi)2^{-j(2H(\xi)+1)}(\log2)(-2j)H'(\xi)\big)(2^{-j}k-t_0)\nonumber\\
&&=K_1(t_0)+\mathcal O(2^{-j(2H(t_0)+1)}j\epsilon_j),
\end{eqnarray}
where $\xi$ is some value in $(\min\{2^{-j}k,t_0\},\max\{2^{-j}k,t_0\})$.

Therefore, by using (\ref{VarV'1}), (\ref{VarV'2}) and the fact that \begin{equation}
\label{identcard}
card(\nu_{t_0,2^j})= 2^{j+1}\epsilon_j+\mathcal O(1),
\end{equation}
 we get
\begin{eqnarray}
\label{VarV'3}
&&\mathbb E\big(V_{X,t_0,j}\big)=c_2(t_0)card(\nu_{t_0,2^j})2^{-j(2H(t_0)+1)}\nonumber\\
&&~~~~+\mathcal O(card(\nu_{t_0,2^j})2^{-j(2H(t_0)+1)}j\epsilon_j)+\mathcal O(2^{-j(2H(t_0)+1)}j^4\epsilon_j)\nonumber\\
&&=2c_2(t_0)2^{-2jH(t_0)}\epsilon_j+\mathcal O\big(2^{-2jH(t_0)}(2^{-j}+j\epsilon_j^2+2^{-j}j^4\epsilon_j)\big).
\end{eqnarray}
Thus  (\ref{VarV'}) follows. Now we are going to show (\ref{VarV}). First observe
$$
Var\big(V_{X,t_0,j}\big)=\sum_{k,k'\in\nu_{t_0,2^j}}Cov\big(d_X(2^{-j},k)^2,d_X(2^{-j},k')^2\big).
$$
Then we recall the fact that if $(Z,Z')$ is a centered Gaussian random vector, then (see e.g. Lemma 5.3.4 in Peng 2011b)
\begin{equation}
\label{identcovGaus}
Cov(Z^2,{Z'}^2)=2\big(Cov(Z,Z')\big)^2.
\end{equation}
It yields
\begin{eqnarray}
\label{VarV1}
&&Var\big(V_{X,t_0,j}\big)=2\sum_{k,k'\in\nu_{t_0,2^j}}
\Big(Cov\big(d_X(2^{-j},k),d_X(2^{-j},k')\big)\Big)^2\nonumber\\
&&=2\sum_{k\in\nu_{t_0,2^j}}\Big(Var\big(d_X(2^{-j},k)\big)\Big)^2\nonumber\\
&&~~+2\sum_{k,k'\in\nu_{t_0,2^j},k\neq k'}\Big(Cov\big(d_X(2^{-j},k),d_X(2^{-j},k')\big)\Big)^2.
\end{eqnarray}
(\ref{VarV1}) together with Lemma \ref{lemma1} implies
\begin{eqnarray}
\label{VarV2}
&&Var\big(V_{X,t_0,j}\big)=2\sum_{k\in\nu_{t_0,2^j}}(c_2(2^{-j}k))^2\big(2^{-j(2H(2^{-j}k) +1)}\big)^2\nonumber\\
&&~~+\mathcal O\Big(\sum_{k\in\nu_{t_0,2^j}}c_2(2^{-j}k)2^{-j(4H(2^{-j}k)+3)}j^4\Big)\nonumber\\
&&~~+
2\sum_{k,k'\in\nu_{t_0,2^j},k\neq k'}(c_1(2^{-j}k,2^{-j}k'))^2\Big(\frac{
2 ^{-j(1+H(2^{-j}k)+H(2^{-j}k'))}}{|k-k'|^{2Q-H(2^{-j}k) -H(2^{-j}k') }}\Big)^2\nonumber\\
&&~~+\mathcal O\Big(\!\!\!\!\sum_{k,k'\in\nu_{t_0,2^j},k\neq k'}\!\!\!\!2 ^{-j(2+2H(2^{-j}k)+2H(2^{-j}k'))}T(k,k',Q,2^{-j},2^{-j})\Big).
\end{eqnarray}
Now we study the second line of (\ref{VarV2}).  By the facts that $\limsup\limits_{j\rightarrow+\infty,k\in\nu_{t_0,2^j}}k2^{-j}=t_0$, $\sup\limits_{j\in\mathbb N,k=0,\ldots,2^j-1}|c_2(2^{-j}k)|<+\infty$, $card(\nu_{t_0,2^j})=\mathcal O(2^j\epsilon_j)$,  we obtain
\begin{equation}
\label{VarV3}
\sum_{k\in\nu_{t_0,2^j}}c_2(2^{-j}k)2^{-j(4H(2^{-j}k)+3)}j^4= \mathcal O(2^{-j(4H(t_0)+2)}j^4\epsilon_j).
\end{equation}
For the fourth line of (\ref{VarV2}), by using Cauchy-Schwarz inequality and Proposition \ref{prop1}, we get
\begin{eqnarray}
\label{VarV4}
&&\sum_{k,k'\in\nu_{t_0,2^j},k\neq k'}2 ^{-j(2+2H(2^{-j}k)+2H(2^{-j}k'))}T(k,k',Q,2^{-j},2^{-j})\\
&&\leq card(\nu_{t_0,2^j})\Big(\!\!\!\!\!\sum_{k,k'\in\nu_{t_0,2^j},k\neq k'}\!\!\!\!\!\!\!\big(2 ^{-j(2+2H(2^{-j}k)+2H(2^{-j}k'))}T(k,k',Q,2^{-j},2^{-j})\big)^2\Big)^{1/2}\nonumber\\
&&\leq c2^{-j(4H(t_0)+3/2)}j\epsilon_j.
\end{eqnarray}
Therefore (\ref{VarV2}) is equivalent to
\begin{eqnarray}
\label{VarV5}
&&Var\big(V_{X,t_0,j}\big)=2\sum_{k\in\nu_{t_0,2^j}}(c_2(2^{-j}k))^2\big(2^{-j(2H(2^{-j}k) +1)}\big)^2\nonumber\\
&&~~~~+
2\sum_{k,k'\in\nu_{t_0,2^j},k\neq k'}(c_1(2^{-j}k,2^{-j}k'))^2\Big(\frac{
2 ^{-j(1+H(2^{-j}k)+H(2^{-j}k'))}}{|k-k'|^{2Q-H(2^{-j}k) -H(2^{-j}k') }}\Big)^2\nonumber\\
&&~~~~+\mathcal O(2^{-j(4H(t_0)+3/2)}j\epsilon_j).
\end{eqnarray}

In order to make an identification of the dominating part of (\ref{VarV5}), for fixed $j,k,k'$ we define $K_2\in C^1([0,1]),~K_3\in C^1([0,1]^2)$ as, for $x,y\in[0,1]$,
\begin{eqnarray*}
&&K_2(x)=(c_2(x))^2\big(2^{-j(2H(x) +1)}\big)^2;\\
&&K_3(x,y)=(c_1(x,y))^2\Big(\frac{
2 ^{-j(1+H(x)+H(y))}}{|k-k'|^{2Q-H(x) -H(y) }}\Big)^2.
\end{eqnarray*}
Then it follows from the Taylor expansion and the fact that $|2^{-j}k-t_0|=\mathcal O(\epsilon_j)$ that
\begin{equation}
\label{VarV6}
K_2(2^{-j}k)=K_2(t_0)+K_2'(\eta)(2^{-j}k-t_0)=K_2(t_0)+\mathcal O(2^{-j(4H(t_0)+2)}j\epsilon_j),
\end{equation}
where $\eta$ is some value in $(\min\{2^{-j}k,t_0\},\max\{2^{-j}k,t_0\})$ and
\begin{eqnarray}
\label{VarV7}
&&K_3(2^{-j}k,2^{-j}k')=K_3(t_0,t_0)+\Big(\frac{\partial K_3(x,2^{-j}k')}{\partial x}\Big|_{x=\eta_1}\Big)(2^{-j}k-t_0)\nonumber\\
&&~~~~+\Big(\frac{\partial K_3(2^{-j}k,y)}{\partial y}\Big|_{y=\eta_2}\Big)(2^{-j}k'-t_0)\nonumber\\
&&=K_3(t_0,t_0)+\mathcal O\bigg(\frac{2^{-j(4H(t_0)+2)}j\epsilon_j}{|k-k'|^{4Q-4H(t_0)}}\bigg),
\end{eqnarray}
where $\eta_1,\eta_2$ are some values satisfying $\eta_1\in(\min\{2^{-j}k,t_0\},\max\{2^{-j}k,t_0\})$ and $\eta_2\in(\min\{2^{-j}k',t_0\},\max\{2^{-j}k',t_0\})$.
Observe that, by a change of variable,
\begin{equation}
\label{VarV21}
\sum_{k,k'\in\nu_{t_0,2^j},k\neq k'}\frac{
1}{|k-k'|^{4Q-4H(t_0)}}= \sum_{|l|\in\nu_{t_0,2^j},l\neq0}\frac{(card(\nu_{t_0,2^j})-|l|)}{|l|^{4Q-4H(t_0)}}.
\end{equation}
 Then (\ref{VarV21}) together with
(\ref{VarV5}), (\ref{VarV6}), (\ref{VarV7}) entails that
\begin{eqnarray}
\label{VarV8}
&&Var\big(V_{X,t_0,j}\big)\nonumber\\
&&=2\Big(c_2(t_0)^2card(\nu_{t_0,2^j})+c_1(t_0,t_0)^2\!\!\!\!\!\sum_{|l|\in\nu_{t_0,2^j},l\neq0}\!\!\!\!\!\frac{(card(\nu_{t_0,2^j})-|l|)}{|l|^{4Q-4H(t_0)}}\Big)2^{-j(4H(t_0)+2)}\nonumber\\
&&~~~~+\mathcal O(2^{-j(4H(t_0)+1)}j\epsilon_j^2)+\mathcal O((2^{-j(4H(t_0)+3/2)}j\epsilon_j)).
\end{eqnarray}
Observe that, since $Q\ge2$,
\begin{eqnarray}
\label{suml}
&&\Big|\sum_{|l|\in\nu_{t_0,2^j},l\neq0}\frac{(card(\nu_{t_0,2^j})-|l|)}{card(\nu_{t_0,2^j})|l|^{4Q-4H(t_0)}}-\sum_{l\in\mathbb Z,l\neq0}|l|^{4H(t_0)-4Q}\Big|\nonumber\\
&&=\mathcal O\Big(\sum_{l=[2^j\epsilon_j]}^{+\infty}l^{4H(t_0)-4Q}+\frac{1}{card(\nu_{t_0,2^j})}\sum_{|l|\in\nu_{t_0,2^j},l\neq0}|l|^{4H(t_0)-4Q+1}\Big).\nonumber\\
\end{eqnarray}
On one hand, we recall that if $f$ is a monotonic decreasing function and the series $\sum\limits_{n=0}^{+\infty}f(n)$ is convergent, then for $N\ge1$,
$$
\int_N^{+\infty}f(x){\,\mathrm{d}} x\leq\sum\limits_{n=N}^{+\infty}f(n)\leq \int_N^{+\infty}f(x){\,\mathrm{d}} x+ f(N).
$$
It yields
\begin{eqnarray}
\label{VarV9}
\sum_{l=[2^j\epsilon_j]}^{+\infty}l^{4H(t_0)-4Q}&=&\mathcal O\Big(\int_{[2^j\epsilon_j]}^{+\infty}x^{4H(t_0)-4Q}{\,\mathrm{d}} x+ (2^j\epsilon_j)^{4H(t_0)-4Q}\Big)\nonumber\\
&=&\mathcal O\big((2^j\epsilon_j)^{4H(t_0)-4Q+1}\big).
\end{eqnarray}
On the other hand, since $4H(t_0)-4Q+1<-1$, then it is easy to see
$$
\frac{1}{card(\nu_{t_0,2^j})}\sum_{|l|\in\nu_{t_0,2^j},l\neq0}|l|^{4H(t_0)-4Q+1}=\mathcal O((2^j\epsilon_j)^{-1}).
$$
Finally,
\begin{equation}
\label{VarV10}
\sum_{|l|\in\nu_{t_0,2^j},l\neq0}\frac{(card(\nu_{t_0,2^j})-|l|)}{|l|^{4Q-4H(t_0)}}=2\Big(\sum_{l\in\mathbb Z,l\neq0}|l|^{4H(t_0)-4Q}\Big)2^j\epsilon_j+\mathcal O(1).
\end{equation}
It follows by (\ref{VarV8}) and (\ref{VarV10}) that
\begin{eqnarray}
\label{VarV11}
&&Var\big(V_{X,t_0,j}\big)=4\Big(c_2(t_0)^2+c_1(t_0,t_0)^2\sum_{l\in\mathbb Z,l\neq0}|l|^{4H(t_0)-4Q}\Big)2^j\epsilon_j2^{-j(4H(t_0)+2)}\nonumber\\
&&~~+\mathcal O(2^{-j(4H(t_0)+2)})+\mathcal O(2^{-j(4H(t_0)+1)}j\epsilon_j^2)+\mathcal O(2^{-j(4H(t_0)+3/2)}j\epsilon_j).
\end{eqnarray}
We thus can conclude
$$
Var\big(V_{X,t_0,j}\big)=c_3(t_0)2^{-j(4H(t_0)+1)}\epsilon_j+\mathcal O(2^{-j(4H(t_0)+1)}(2^{-j}+j\epsilon_j^2+2^{-j/2}j\epsilon_j)\big),
$$
where $c_3(t_0)=4(c_2(t_0)^2+c_1(t_0,t_0)^2\sum_{l\in\mathbb Z,l\neq0}|l|^{4H(t_0)-4Q})>0$ is a constant only depending on $t_0$. Proposition \ref{VarVX} has been proven. $\square$

\end{document}